\documentclass[11pt]{amsart}
\usepackage[utf8]{inputenc}
\usepackage{amsmath,amsxtra,amssymb,latexsym,amscd,amsthm,bbm}
\usepackage[colorlinks,citecolor=blue,linkcolor=blue,pagebackref]{hyperref}
\usepackage{url}
\usepackage{amsaddr} 
\usepackage{xcolor}
\usepackage{tikz}
\usepackage{pgfplots}
\usepackage{enumitem}

\usepackage{graphicx}
\usepackage{caption}
\usepackage{subcaption}

\newtheorem{theorem}{Theorem}[section]
 
\newtheorem{lemma}[theorem]{Lemma}
\newtheorem{definition}[theorem]{Definition}
\newtheorem{remark}[theorem]{Remark}
\newtheorem{assumption}[theorem]{Assumptions}

\definecolor{darkgreen}{rgb}{0,0.5,0}

\newcommand{\R}{\mathbb{R}}
\newcommand{\ud}{\mathrm{d}}

\newcommand{\gep}{\varepsilon}

\title[Particle trajectories of scalar conservation laws and BIPs]{Stability of particle trajectories of scalar conservation laws and applications in Bayesian inverse problems}
\author[Masoumeh Dashti and Duc-Lam Duong]{Masoumeh Dashti$^{\dag}$ and Duc-Lam Duong$^{\ddag *}$} 
\address{$^{\dag}$Department of Mathematics\\
School of Mathematical and Physical Sciences, University of Sussex\\
Brighton, UK\\
e-mail: m.dashti@sussex.ac.uk
\\ and \\
$^{\ddag}$Department of Computational Engineering\\
School of Engineering Science, LUT University\\
Lappeenranta, Finland\\
e-mail: duc-lam.duong@lut.fi}
\thanks{$^{*}$Corresponding author.}

\begin{document}

\begin{abstract}
We consider the scalar conservation law in one space dimension with a genuinely nonlinear flux. We assume that an appropriate velocity function depending on the entropy solution of the conservation law is given for the comprising particles, and study their corresponding trajectories under the flow. The differential equation that each of these trajectories satisfies depends on the entropy solution of the conservation law which is typically discontinuous in both time and space variables. 
The existence and uniqueness of these trajectories are guaranteed by the Filippov theory of differential equations. We show that such a Filippov solution is compatible with the front tracking and vanishing viscosity approximations in the sense that the approximate trajectories given by either of these methods converge uniformly to the trajectories corresponding to the entropy solution of the scalar conservation law. For certain classes of flux functions, illustrated by traffic flow, in our main result, we prove the H\"older continuity of the particle trajectories with respect to the initial field or the flux function. 
 We then consider the inverse problem of recovering the initial field or the flux function of the scalar conservation law from discrete pointwise measurements of the particle trajectories. We show that the above continuity properties translate to the stability of the Bayesian regularised solutions of these inverse problems with respect to appropriate approximations of the forward map. We also discuss the limitations of the situation where the same inverse problems are considered with pointwise observations made from the entropy solution itself.\\

    
\end{abstract}

\maketitle

\tableofcontents

\section{Introduction}
\subsection{Scalar conservation laws}
We consider the scalar conservation law in one space dimension
\begin{equation} \label{scl}
\partial_t v(x,t) + \partial_x f(v(x,t)) = 0,  \quad x\in \mathbb{R},\; t>0,
\end{equation}
with an initial condition
\begin{equation}\label{ini}
v(x,0) = v_0(x),  \quad x\in \mathbb{R}.
\end{equation}
Here \(v:\mathbb{R} \times [0,\infty) \to \mathbb{R}\) denotes the density or concentration of some conserved physical quantity, and \(f: \mathbb{R} \to \mathbb{R}\) represents its flux. 
 An outstanding example of such conservation laws is the LWR model of traffic flow on a highway, initiated by Lighthill-Whitham (\cite{LW55}) and Richards (\cite{Ric56}). In that model, $v(x,t)$ (which is often denoted by $\rho(x,t)$ in this context) represents the density of the vehicles at location $x$ and time $t$, measured by the number of vehicles per unit length. The quantity $f(v(x,t))$ is the flux of vehicles across the point $x$ on the road at time $t$. 

It is well-known that conservation law \eqref{scl}-\eqref{ini} exhibits \emph{shocks}, even if the flux $f$ and initial data $v_0$ are smooth. Due to this, solutions must be sought in the space of discontinuous functions, and weak solutions must be used. Weak solutions of conservation laws are generally not unique unless additional constraints, referred to as \emph{entropy conditions}, are imposed. The resulting solution is called an \emph{entropy solution}. In the case of scalar conservation laws, $v$ is an entropy solution to \eqref{scl} if it satisfies the inequality \begin{equation*}
    \partial_t \eta(v(x,t)) + \partial_x q(v(x,t)) \le 0,
\end{equation*}
in the sense of distribution for all entropy-entropy flux pairs $(\eta,q)$ where $\eta$ is convex and $q'(v) = \eta'(v)f'(v)$. Much research has been devoted to studying the existence, uniqueness, and regularity of the entropy solution to \eqref{scl}-\eqref{ini}, given appropriate flux and initial data. For a comprehensive treatment of the subject, refer to the book by Dafermos \cite{Daf16}, as well as the monographs by Bressan \cite{Bre00} and by Holden and Risebro \cite{HR15}.

Motivated by the traffic flow model, in this work, we assume that the physical system that \eqref{scl}-\eqref{ini} describes is comprised of particles moving with the flow. We denote by $w$ the velocity of the flow which we suppose to be a function of only $v$. The trajectory $z(t)$ of a particle starting from some point $x_0$ is then defined by the ordinary differential equation
\begin{equation}\label{ode}
     \frac{dz(t)}{dt} = w(v(z(t),t)),
\end{equation}
subject to the initial condition $z(t_0) = x_0$, $t_0>0$. We explain below, in the second part of this introduction, that one motivation for studying such trajectories is their application in the inverse problems in the context of scalar conservation laws.

Since $v$ is typically a discontinuous function of both $z$ and $t$, the standard Cauchy-Lipschitz theory for ODEs does not apply to \eqref{ode}. In this work, we employ the Filippov theory \cite{Fil88} which defines an absolutely continuous function $z$ to be a solution to \eqref{ode} if $z$ satisfies \eqref{ode}, viewed as a differential inclusion, almost everywhere (see Definition \ref{d:FilS}). 
We are interested in the stability properties of $z$, the Filippov solution of \eqref{ode}, with respect to the initial field $v_0$ or the flux function $f$ in the equation \eqref{scl}. The existing theory (using Oleinik's decay estimate, see, for instance, \cite{Leg11}) proves the existence and uniqueness of $z$, but the method cannot be used to prove the stability of $z$ with respect to changes in the velocity field or the flux function. We establish these stability properties by combining the Filippov theory, the front tracking method, and some structural properties of solutions of conservation laws (see Section \ref{subsec:Lagrangian}, Theorem \ref{thm:trajectory_cont}). 

We also consider the approximate trajectories arising from the method of vanishing viscosity. Let $v^\epsilon$ be the solution of the viscous version of \eqref{scl}, that is
\begin{equation*} 
\partial_t v^\epsilon(x,t) + \partial_x f(v^\epsilon(x,t)) = \epsilon \partial_{xx} v^\epsilon(x,t).
\end{equation*}
Consider the particle trajectory $z^\epsilon$ starting at $z^\epsilon(0) = x_0$ that solves
\begin{equation*}
    \dot z^\epsilon(t) = w(v^\epsilon(z^\epsilon(t),t)).
\end{equation*}
We show that as $\epsilon \to 0$, the trajectory $z^\epsilon$ converges in $L^\infty$  to the Filippov solution $z$ of \eqref{ode}. To prove this result, we make the assumption that the trajectory $z(t)$, even though it may cross the shock curves of $v$, never lies on any shock curve for a positive period of time (Assumption \ref{Ass:shockspeed}). We then verify that this assumption is at least satisfied for the case of traffic flow provided that the initial field is strictly positive.

We then study in detail the traffic flow model, the example that motivates our work. Denoting the vehicle density by $\rho$ and the vehicle speed by $w$, the equation describing the traffic flow \cite{LW55, Ric56} is written as
\begin{equation} \label{eq:traffic}
    \partial_t\rho + \partial_x (\rho w(\rho)) = 0,\qquad\rho(\cdot,0)=\rho_0.
\end{equation}
In this case, we are able to obtain stronger stability results. We show $1/2$-H\"older continuity of the mapping $\rho_0\mapsto z(\cdot): L^1\cap BV\to L^\infty$ (see Theorem \ref{thm:rate_traffic}) where $BV$ denotes the space of functions of bounded variations, provided that the function $w$ is strictly decreasing. A similar stability estimate with respect to changes in flux function is also proved (Theorem \ref{thm:rate_flux}). 
In both cases, the main work is in estimating the error incurred when the two trajectories we are comparing pass through shocks of relatively large size. We estimate the error for a single shock and show that, for an initial field of bounded variation, the sum of the errors after passing all the shocks in a finite time interval remains small.

\subsection{Bayesian inverse problems in scalar conservation laws}
In the second part of the paper, we study the inverse problems of recovering the initial field $v_0$ (with $f$ given) or the flux function $f$ (with $v_0$ given) from observations of (a function of) the solution $v$. These kinds of inverse problems have many applications, depending on how we interpret $v$ in the model \eqref{scl}. In the case of traffic flow, they correspond to the problems of determining the upstream vehicle density or the flux given finite observations of (known functions of) the density field at later times.
In the language of mathematics, these inverse problems can be written as 
\begin{equation*}
y = G(u),
\end{equation*}
where \(u\) denotes the \emph{unknown} (that is, $u$ is either $v_0$ or $f$), \(y\) is the \emph{observed data} and \(G\) is the \emph{observation map} which is defined via the solution of the forward problem. We consider the observation map to be 
\begin{equation}\label{LagDa}
    G(u) = \{z(t_j)\}_{j\in J},
\end{equation}
for some finite index sets $J$. 
In other words, we gather the observations by tracking the trajectory of a particle moving along the flow at discrete times. The treatment of the case where the observations are made from more than one trajectory would be similar.

We note that, in general, inverse problems are ill-posed, meaning solutions may not exist, may not be unique or may depend sensitively on data. In hyperbolic conservation laws, the situation is particularly complicated since the physically relevant solutions are often irreversible. 
This irreversibility property, induced by the entropy condition (see \cite{Daf16}), renders severe difficulties in inversion: 
if $v(x,t)$ is an entropy solution, $v(-x,-t)$ is no longer an entropy solution unless $v(x,t)$ is a classical one (and has in particular no shocks). Finding the right techniques to tackle inverse problems in hyperbolic conservation laws is therefore challenging.

Nevertheless, by appropriate regularisation one is able to find some estimation of the missing information. 
We employ a Bayesian approach to regularisation which, as the solution, provides a probability distribution on the unknown called the posterior. 
The Bayesian inverse problems for unknown functions have been studied extensively in the last decade, in particular, for nonlinear models involving PDEs; see \cite{Stu10, DS16}, the early paper \cite{Fra70}, and for a more applied and computational overview, \cite{KS05}. 
In the Bayesian framework, the data and the unknown are treated as random variables, and the regularisation enters the framework in the form of a given prior probability distribution on the unknown. The posterior may then be derived through Bayes' theorem and depends on the prior and the data and also on the forward map through the observation operator. 
The continuity properties of the particle trajectories mentioned above suggest that if we consider the observation map as \eqref{LagDa}
then the collected data is stable in some sense to be made precise later on, providing some regularity structure for the observation map.  Thanks to the approximation theory of Bayesian inverse problems (see Section \ref{sec:Bay_recovery}), the approximate posterior can be shown to be continuous in appropriate metrics, giving the well-posedness for the solutions of our inverse problems. 
In the case of the traffic flow, the strong stability properties of the trajectory $z$ (Theorem \ref{thm:rate_traffic} and \ref{thm:rate_flux}) will be of great value for our inverse problems since they translate to the rate of convergence for the approximations of the corresponding posteriors.

In traffic flow applications for example, a different set of data may be available. One may have discrete measurements of the entropy solution $v$ itself, that is
\begin{equation*}
    G(u) = \{v(x_i,t_j)\}_{i\in I, j\in J},
\end{equation*}
for finite sets $I, J$. In the context of conservation laws, due to the discontinuities in $v$, the observation operator lacks desirable regularity properties for the approximation theory of the Bayesian approach employed in \cite{Stu10} to work. Nevertheless, we still have a well-posedness result for the Bayesian inverse problems, thanks to the measurability of the forward map. See Section \ref{sec:app_Euler} for a more detailed discussion. 

We note that, in many situations, making observations by tracking particle trajectories is probably more practical and more economical than measuring the flow field itself. Consider for example the traffic flow passing through a tunnel, where measuring the density of the cars inside the tunnel might not be easy, one may instead track the position of a marked car over time. In practice, this can be done easily via a GPS device mounted on the car. If necessary, at the same time one can track more cars to have a more accurate picture.

\subsection{Our contributions}
Let us summarise here the main contributions of this paper, which are twofold. 

\begin{itemize}
    \item The first main contribution of our paper is the study of the particle trajectories \eqref{ode}, in connection with the entropy solution to \eqref{scl}-\eqref{ini}. Related research in this direction traced back to \cite{BS98} where the authors considered \eqref{ode} with the right-hand side connected to a $2\times 2$ system of conservation law (see also \cite{bressan1988unique} for an earlier work on discontinuous ODEs). Other works where the equation \eqref{ode} is motivated from a traffic flow model were considered in \cite{CM03, Mar04}, and later on in \cite{DMG14}. We note, however, that \cite{DMG14} only considers solutions in the sense of Carathéodory, which is somewhat more restricted than those considered in this paper, while \cite{CM03} focuses on the traffic flow model.
    In this work, we are able to provide strong stability results that, to our knowledge, have not been previously studied in the literature.

    \item The second main contribution goes towards the formulation and well-posedness establishment of the Bayesian inverse problems for scalar conservation laws. Due to the nature of shockwaves, inverse problems for conservation laws are challenging and works in this direction are rather limited in the literature. Of recent contributions toward understanding these inverse problems, we mention the identification problem of the (possibly discontinuous) flux function considered in \cite{HPR14} and of the initial data in \cite{colombo2020initial} and \cite{liard2021initial}. The well-posedness of the Bayesian inverse problems for hyperbolic conservation laws is also recently considered in \cite{mishra2021well} (however with a totally different forward map) and somewhat in \cite{lanthaler2022bayesian} where the focus is on the data assimilation problem.
\end{itemize}

To conclude the introduction, we note here that the notion of Filippov solutions to differential equations with discontinuous right-hand side was already employed by Dafermos in \cite{Daf77} to build the theory of generalized characteristics for hyperbolic conservation laws and has been an efficient method for studying the regularity of solutions. It is worth noting that the trajectories that we consider in this paper differ from the generalized characteristics considered by Dafermos. The speed of the generalized characteristics is either the classical characteristics speed or shock speed, while the particle speed considered here, given on the right-hand side of \eqref{ode}, is the speed of the flow itself. 

\subsection{Organisation of the paper} The paper is organised as follows. In Section \ref{sec:SCL_Lagrangian}, after recalling a few basic notions and properties of the solutions, we show the continuity of particle trajectories with respect to appropriate approximations of the solutions of scalar conservation laws. In Section \ref{sec:stab_Lagrangian}, in the case of traffic flow, we prove H\"older stability of vehicle trajectories with respect to the initial field and flux function.
The Bayesian inverse problems for the initial field or flux function, given discrete noisy observations of a particle trajectory, is considered in Section \ref{sec:Bay_recovery}.
In the last part of Section \ref{sec:Bay_recovery}, we will discuss the case where the data comes from pointwise measurements of the entropy solution itself. 

    

\section{Scalar conservation laws and the particle trajectories}\label{sec:SCL_Lagrangian} 

%
\subsection{Basic notions and properties}
In this section, we recall some fundamental notions and properties of scalar conservation laws. Consider the Cauchy problem
\begin{equation} \label{scl-2}
\partial_t v(x,t) + \partial_x f(v(x,t)) = 0,  \quad x\in \mathbb{R}, t>0,
\end{equation}
\begin{equation}\label{ini-2}
v(x,0) = v_0(x),  \quad x\in \mathbb{R},
\end{equation}
where $f$ is assumed to be at least locally Lipschitz continuous and $v_0$ is a bounded measurable function. The development of shockwaves in general leads to the consideration of weak solutions.
\begin{definition}[Weak solution]\label{def:weak_sol}
A function \(v \in L^\infty(\mathbb{R} \times (0,\infty))\) is called a weak solution for \eqref{scl-2}-\eqref{ini-2} if
\begin{equation*}
\int_0^\infty \int_{\mathbb{R}} v(x,t) \phi_t + f(v(x,t)) \phi_x dx dt = 0
\end{equation*}
holds for any test function \(\phi \in C_0^\infty(\mathbb{R} \times (0,\infty))\).
\end{definition}
Weak solutions to \eqref{scl-2}-\eqref{ini-2} are, therefore, allowed to have discontinuities. Nevertheless, not all discontinuities are permitted except the ones that satisfy the following jump condition.\\
{\bf Rankine-Hugoniot condition.} On every discontinuity curve $\alpha = \alpha(t)$, we have
\begin{equation}\label{Rankine-Hugoniot}
    s: = \alpha'(t) = \frac{f(v_l) - f(v_r)}{v_l - v_r},
\end{equation}
where $v_l$ and $v_r$ denote the limits from the left and right of the discontinuity curve. The term $s$ is often referred to as \emph{shock speed}.

The notion of weak solution, however, is too weak to ensure uniqueness. To single out the physically relevant solution, one needs to fill in extra information.
\begin{definition}[Entropy solution]\label{def:entropy}
For every convex function \(\eta\) (which is called an entropy), we define the entropy flux \(q\) by \(q'(u)=\eta'(u)f'(u)\). A function \(v \in L^\infty(\mathbb{R} \times (0,\infty))\) is called an entropy solution for \eqref{scl-2}-\eqref{ini-2} if the inequality
\begin{equation*}
\int_0^\infty \int_{\mathbb{R}} \eta(v(x,t)) \phi_t + q(v(x,t)) \phi_x dx dt \ge 0
\end{equation*}
holds for any convex entropy/entropy flux pair \((\eta, q)\) and any non-negative test function \(\phi \in C_0^\infty(\mathbb{R} \times (0,\infty))\).
\end{definition}

With the entropy condition being added, solutions to \eqref{scl}-\eqref{ini} become globally well-posed thanks to the following classical result.
\begin{theorem}[Kruzkov \cite{Kru70}] \label{thm:Kruzkov} 
For every \(v_0 \in L^\infty(\mathbb{R})\), there exists a unique entropy solution $v$ to \eqref{scl}-\eqref{ini} in $C([0,\infty); L^1_{\mathrm{loc}}(\mathbb{R}))$ that satisfies, for every $t>0$,
\begin{equation}\label{eq:infbnd}
    \|v(\cdot,t)\|_{L^\infty(\mathbb{R})} \le \|v_0\|_{L^\infty(\mathbb{R})}.
\end{equation}
Moreover, if $\bar{v}$ is the entropy solution to \eqref{scl} corresponding to the initial data $\bar{v}_0$ with $ v_0, \bar{v}_0 \in L^1(\mathbb{R}) \cap L^\infty(\mathbb{R})$, then 
\begin{equation}\label{eq-L1stab1}
\|v(\cdot,t) - \bar{v}(\cdot,t)\|_{L^1(\mathbb{R})} \le \|v_0-\bar{v}_0\|_{L^1(\mathbb{R})},
\end{equation}
for all $t>0$.
\end{theorem}

\subsubsection{Vanishing viscosity approximation}
The entropy solution defined via Definition \ref{def:entropy} agrees with the viscosity solution obtained via the vanishing viscosity method (and hence, it gives a physically meaningful solution, see a discussion in \cite[Chapter 4]{Daf16}). That is, if we add a (small) diffusion term to the right-hand side of \eqref{scl} and consider a solution to the new equation, then at the limit when the diffusion coefficient vanishes, one gets the entropy solution to \eqref{scl}. Indeed, consider the parabolic equation
\begin{equation} \label{eq_viscosity}
\partial_t v(x,t) + \partial_x f(v(x,t)) = \epsilon\partial_{xx} v(x,t),
\end{equation}
with some $\epsilon > 0$. Assume that \(v^\epsilon\) is a smooth solution of \eqref{eq_viscosity}-\eqref{ini} (which always exists and is unique by the regularity of parabolic equations). Multiply both sides of \eqref{eq_viscosity} by $\eta'(v^\epsilon(x,t))$ to get 
\begin{equation*} 
\partial_t \eta(v^\epsilon) + \partial_x q(v^\epsilon) = \epsilon\partial_{xx} \eta (v^\epsilon) - \epsilon \eta''(v^\epsilon)|\partial_x v^\epsilon|^2 \le \epsilon\partial_{xx} \eta (v^\epsilon).
\end{equation*}
Thus if we let \(\epsilon \to 0\) and assume that,
\[v^\epsilon \to v \quad \text{a.e. boundedly,}\]
then \(v\) satisfies Definition \ref{def:entropy} and is an entropy solution to the scalar conservation law \eqref{scl}-\eqref{ini}.

\subsubsection{Front tracking approximation}
\label{subsec:FTA}
The front tracking method, introduced by Dafermos \cite{Daf72} and developed by DiPerna \cite{DiP76}, Bressan \cite{Bre92} and Risebro \cite{Ris93}, is a powerful tool in the existence theory of entropy solutions in both the scalar case and systems of hyperbolic conservation laws. 
The idea is to approximate the initial function $v_0$ (of bounded variation) by a step function $v_0^N$ and the flux $f$ by a piecewise linear function $f^N$. The approximated solution $v^N$ is then given by solving a (finite) set of so-called Riemann problems, each problem is associated with a point of discontinuity in initial data, given as

\begin{equation}
\label{Riemann_ini}
 v(x,0) =  \begin{cases}
v_l & \text{ if } x < x_j \\ 
v_r & \text{ if } x > x_j,
\end{cases}
\end{equation}
where $x_j$ is a point of discontinuity and $v_l$ and $v_r$ denote the values of $v_{0,N}$ at the left and the right limits $v_{0,N}(x_j-), v_{0,N}(x_j+)$.
One claims that Riemann problems obey a maximum principle, meaning the solution to \eqref{scl_2} with initial data \eqref{Riemann_ini} remains in between $v_l$ and $v_r$, and that the solution of \eqref{scl_2} will take values in the set $\{v_{j,N}\} \cup \{\text{break points of } f^N_\smile \text{ or } f^N_\frown\}$, where $f_\smile$ (or $f_\frown$) denotes the convex envelope (or concave envelope) of $f$. This solution is defined up to some time $t = t_1$ where two or more jump discontinuities (coming from nearby Riemann problems) collide, forming new Riemann problems. The above procedure continues for new Riemann problems and the solution is prolonged up to some new collision time $t=t_2$, and so on. Luckily enough, this process does not go on forever thanks to the fact that the number of interactions is finite, see for example \cite{Bre00} for a proof. The intuition behind this fact can be understood roughly like this, each time a new collision forms, two or more discontinuities collapse to produce a single discontinuity, thus the wave pattern is simplified since the number of jump discontinuities is decreasing over time. In particular, if the initial condition $v_0$ is a non-negative differentiable function with compact support and the flux function is smooth and uniformly convex, then after a certain time, all shocks will finally be merged and continue as a single shock (see \cite{KT05}, also \cite{Whi75}). 

The method of front tracking approximations provides an alternative approach for proving the existence and uniqueness of entropy solutions of \eqref{scl}-\eqref{ini}. Moreover, one can use this method to derive the stability property of the entropy solution with respect to the flux function, as the following result shows. Hereafter we denote 
\begin{equation}
\label{BV}
    BV(\mathbb{R}) := \{u: \mathbb{R} \to \mathbb{R}, {\rm TV}(u) < \infty\},
\end{equation}
the space of bounded variation functions, with ${\rm TV}(u)$ being the total variation of the function $u$; and
\begin{equation}
\label{Lip-norm}
    \|f\|_{\rm Lip} := \sup_{u\neq v} \left| \frac{f(u) - f(v)}{u-v}\right|,
\end{equation}
the Lipschitz constant for $f$.
\begin{theorem} \label{thm:Lucier}
Let \(v_0 \in L^1 \cap BV(\mathbb{R}) \) and $v, \bar{v}$ be the entropy solutions to \eqref{scl}-\eqref{ini} with respect to locally Lipschitz continuous flux functions $f, g$ (respectively), then there exists a constant \(C=C(u_0)\) such that 
\begin{equation} \label{eq:L1stab2}
\|v(\cdot,t) - \bar{v}(\cdot,t)\|_{L^1(\mathbb{R})} \le Ct\|f-g\|_{\rm Lip},
\end{equation}
for every \(t>0\).
\end{theorem}
A proof of this theorem using front tracking method can be found in \cite[Chapter 2]{HR15}. The assumption $v_0 \in BV(\mathbb{R})$ is needed in establishing \emph{a priori} a bound on ${\rm TV}(v(\cdot,t))$ which is an essential part in the proof of existence (typically via Helly's theorem). The uniqueness property \eqref{eq:L1stab2}, was first obtained in \cite{Luc86} with a slightly different approach.

\subsection{Particle trajectories}
\label{subsec:Lagrangian}
Let $v$ be the unique entropy solution of the scalar conservation law
\begin{equation}\label{scl_2}
    \partial_t v + \partial_x f(v) = 0,
\end{equation}
with initial data $v_0 = v(\cdot,t_0) \in [-M, M]$ and locally Lipschitz continuous flux function $f$ which is assumed to be genuinely nonlinear, that is, $f$ is either strictly convex or strictly concave on $[-M, M]$. 

 We assume that the physical system that equation (\ref{scl_2}) describes,  is comprised of particles moving with the flow field. 
Let us consider a single particle starting at $t_0$ from a point $x_0$. Denote by $z(t)$ its position at time $t$ (and, by abuse of language, we will also call the particle itself by $z$). 
Then $z$ satisfies the following equation
\begin{align}\label{e:LagS}
    \dot{z}(t) = w(v(z(t),t)),\quad  z(t_0) = x_0
\end{align}
where $w$ satisfies the following assumption.
\begin{assumption}\label{as:genw}
The velocity field $w$ is a bounded and continuous function of $v$ which is non-increasing if $f$ is strictly concave, and non-decreasing if $f$ is strictly convex.
\end{assumption}
This assumption is satisfied, in particular, for the traffic flow model studied in Section \ref{sec:stab_Lagrangian}.
For that model, $v$ represents the density
of cars, where a typical example of the velocity of the car flow $w$ is given by $w(v) = w_{max}(1-v/v_{max})$ (see the discussions in Section \ref{subsec:traffic}).

We note that, since $v$ is a solution of the conservation law \eqref{scl_2}, the function $\alpha(x,t):= w(v(x,t))$ on the right-hand side of the ODE \eqref{e:LagS} might be discontinuous in both $x$ and $t$, hence the Cauchy–Lipschitz theory for ODEs does not apply. Therefore more attention should be paid to defining the trajectory $z$ as a solution of \eqref{e:LagS} in an appropriate way. Here we use the following definition by Filippov \cite{Fil88}.
\begin{definition}\label{d:FilS}
A function $z: [t_0, T] \to \mathbb{R}$ is called a solution to \eqref{e:LagS} in the sense of Filippov if it is absolutely continuous on $[t_0, T]$ and it satisfies the differential inclusion
\begin{equation}\label{Filippov_incl}
    \dot{z}(t) \in [w(v(z(t)\pm,t)), w(v(z(t)\mp,t))],~~~\mbox{ for almost every $t\in [t_0, T]$}
\end{equation}
for strictly convex and strictly concave flux function $f$ respectively, where $v(x\pm,t)$ denote the one-sided limits of $v$ at $x$.
\end{definition}

The traces $v(x\pm,t)$ exist thanks to a classical result that $v(\cdot,t)\in BV_{loc}(\mathbb{R})$ for all $t>0$ (even if $v_0$ is merely in $L^\infty$, see \cite{Daf16}, Chapter XI). Moreover, it follows from the Lax entropy condition that
$$v(x+,t) \le v(x-,t),$$ 
for almost all $t > 0$ and all $x\in \mathbb{R}$, when $f$ is strictly convex and with the reverse inequality when $f$ is strictly concave. Hence the right-hand side of \eqref{Filippov_incl} makes sense. 
The following theorem is an easy application of Filippov's theory and ensures the existence and uniqueness of $z(t)$. The result is already obtained in \cite{CM03} for concave flux functions and, with a slightly different setting, in \cite{Leg11}. 

\begin{theorem} \label{thm:uniq_trajectory}
Let v be a unique entropy solution of \eqref{scl}-\eqref{ini} with $f$ being strictly convex or strictly concave. Then for $(x_0,t_0) \in \mathbb{R} \times (0,\infty)$, there exists a unique absolutely continuous function $z: [t_0,\infty) \to \mathbb{R}$ satisfying \eqref{e:LagS} in the sense of Filippov.
\end{theorem}

\begin{remark}
By using the technique in the proof of Theorem \ref{thm:uniq_trajectory},  one can prove the following stability-like estimate 
\begin{equation}\label{eq:stab_ini_posit}
    |x(t) - y(t)|^2 \le |x_0 - y_0|^2\left(\frac{t}{t_0}\right)^C,
\end{equation}
for $x, y$ being the Filippov solution of \eqref{e:LagS} with respect to initial positions $x_0, y_0$, respectively, where $C$ depends only on $f$ and the Lipschitz constant of $w$. 
However, this proof cannot be used to derive the continuity or stability of the particle trajectories with respect to changes in the velocity field $v$ itself as a result of the changes in the initial field which we study here. This is because now the two trajectories do not solve the same ODE anymore. 
\end{remark} 

In the rest of this section, for convenience, we will assume that $f$ is strictly convex, with notice that the analysis applies to the case of concave flux functions as well. The velocity field $w(v)$ is then assumed to be a non-decreasing function of $v$. 

We now turn to investigate the continuity properties of the particle trajectories, defined as the Filippov solutions of (\ref{e:LagS}), with respect to changes in the solution field $v$. These changes may be a result of perturbations in the initial field or from any kind of approximations of the forward model incurred in a computational process. In particular, we consider two of the most popular approximations used in scalar conservation laws, the front tracking approximation and the vanishing viscosity approximation. We prove that the approximations of solution of (\ref{e:LagS}) arising from a small perturbation of the entropy solution of (\ref{scl_2}) (including the front tracking approximations) converge to the unique Filippov solution of \eqref{e:LagS}. We prove a similar stability property when the approximation is a result of the vanishing viscosity method, however, with some restriction on the shock speed.

In this context, it is natural to consider some appropriate notion of approximate solutions. Following Filippov, for some set $K$, we denote by $K^\delta$ the following set 
\[ K^\delta := \left\{x: \inf_{y\in K}|y-x| \le \delta \right\}, \]
that is, a closed $\delta$-neighbourhood of $K$. 
\begin{definition}
We call $y(t)$ a \emph{$\delta$-solution} of the inclusion \eqref{Filippov_incl} if $y(t)$ is absolutely continuous and we have, almost everywhere,
\begin{equation}
    \dot{y}(t) \in V(y(t)^\delta,t)^\delta,  
\end{equation}
where $V(y,t) := [w(v(y+,t)), w(v(y-,t))]$.
\end{definition}

The following lemma will be of later use. For the proof we refer to Filippov \cite[Chapter 2, \S 7]{Fil88}.
\begin{lemma}\label{lem:Fil-approx}
Let $x_k(t)$ be a uniformly convergent sequence of $\delta_k-$solutions of the inclusion \eqref{Filippov_incl} with $\delta_k \to 0$ as $k\to \infty$. Then the limit $x(t) = \lim_{k\to \infty} x_k(t)$ is also a solution of this inclusion.
\end{lemma}

We now consider a sequence of exact solutions (or front tracking approximations) $v^N$ that converges to the solution $v$ in $L^1$. The Filippov solution to \eqref{e:LagS} with $v$ replaced by $v^N$ is denoted by $z^N$, that is, $z^N$ satisfies
\begin{equation}\label{eq:Lagrangian_appx}
    \dot{z}^N(t) \in [w(v^N(z^N(t)+,t)), w(v^N(z^N(t)-,t))],
\end{equation}
for almost every $t$. We assume that the starting position of $z^N$ and the starting position of $z$ are the same,  $z^N(t_0) = z(t_0) = x_0$. The following result establishes the uniform convergence for $z^N$.

\begin{theorem} \label{thm:trajectory_cont}
Let $v^N$ be a sequence of exact solutions (or front tracking approximations) of bounded variation converging in $L^1$ to the entropy solution $v$ of \eqref{scl_2}. Let $z^N$ be defined as \eqref{eq:Lagrangian_appx}. Then $z^N$ converges to $z$ uniformly on $[t_0, T]$ for every $T>t_0>0$, as $N\to \infty$, where $z$ is the Filippov solution of \eqref{e:LagS}.
\end{theorem}

\begin{proof}[Proof of Theorem \ref{thm:trajectory_cont}] We first consider the case where $\{v^N\}$ is a front tracking approximation of $v$. 
Fix $T>0$. According to \eqref{e:LagS} and \eqref{eq:infbnd}, and since $w$ is bounded and continuous by Assumption \ref{as:genw}, there exists a constant $C$ such that, for any $t\in [t_0,T]$,
\begin{equation*}
    |\dot{z}^N(t)| \le C.
\end{equation*}
This implies that $\{z^N\}$ is uniformly bounded and also
\begin{equation*}
    |z^N(t) - z^N(s)| \le C|t-s|,
\end{equation*}
which means $\{z^N\}$ is an equicontinuous sequence. By Arzela-Ascoli theorem, there exists a subsequence, still denoted by $z^{N}$, such that 
\begin{equation} \label{uniform_conv_2}
    z^{N}(\cdot) \to z(\cdot) \quad \text{ in } C^0([t_0,T]),
\end{equation}
for some Lipschitz continuous function $z$. We claim that $z$ solves \eqref{e:LagS} in the sense of Filippov, that is, a.e.-$t$,
\begin{equation*}
    \dot{z}(t) \in V(z(t),t) = [w(v(z(t)+,t)), w(v(z(t)-,t))].
\end{equation*}
Indeed, from the definition of $z^N$ we have, a.e.-$t$,
\begin{equation*}
    \dot{z}^N(t) \in V^N(z^N(t),t) := [w(v^N(z^N(t)+,t)), w(v^N(z^N(t)-,t))].
\end{equation*}
We now proceed as in \cite[Section 4.3.2]{DMG14} (see also \cite{BL99}) to prove that, a.e. $t \in [t_0,T]$,
\begin{equation}\label{eq:51} 
    v^N(z^N(t)+,t) \to v^+(t) := v(z(t)+,t), \quad \text{as } N\to\infty.
\end{equation}
Indeed, by extracting a further subsequence if needed, $v^N$ converges a.e. to $v$, there exists a sequence $\bar{z}^N \ge z^N(t)$ such that $\bar{z}^N \to z(t)$ and $v^N(\bar{z}^N,t) \to v^+(t)$. 

For a.e. $t$, if $(z(t),t)$ is a point of continuity of $v$ then, for any fixed $\varepsilon_0 > 0$, there exists $\delta > 0$ such that $TV(v(\cdot,t): (z(t)-\delta, z(t)+\delta)) \le \varepsilon_0$. Then
\begin{equation*}
    TV(v^N(\cdot,t): (z(t)-\delta, z(t)+\delta)) \le 2\varepsilon_0,
\end{equation*}
for large enough $N$, by weak convergence of measures (see \cite[Lemma 15]{BL99}). Therefore
\begin{equation*}
    |v^N(z^N(t)+,t) - v^+(t)| \le |v^N(z^N(t)+,t) - v^N(\bar{z}^N,t)| + |v^N(\bar{z}^N,t) - v^+(t)| \le 3\varepsilon_0,
\end{equation*}
for large enough $N$. Since $\varepsilon_0$ can be chosen to be arbitrary small, it implies \eqref{eq:51}.

If $(z(t),t)$ is a point of discontinuity of $v$ with $|v(z(t)+,t) - v(z(t)-,t)| \ge \varepsilon_0$, then also
\begin{equation*}
    |v^N(z^N(t)+,t) - v^N(z^N(t)-,t)| \ge \frac{\varepsilon_0}{2},
\end{equation*}
for $n$ large enough. We will prove that, for each $\varepsilon > 0$, there exists $\delta > 0$ such that for all $n$ large enough we get
\begin{equation}\label{eq:52}
    |v^N(x,s) - v^N(z^N(t)+,t)| < \varepsilon, \;\; \text{for } |s-t|\le \delta,\; |x-z(t)|\le \delta, \; x>z^N(s).
\end{equation}
Indeed, if \eqref{eq:52} does not hold, there will be $\varepsilon > 0$ and sequences $t_N \to t$, $\delta_N \to 0$ such that
\begin{equation*}
    TV(v^N(\cdot, t_N): (z^N(t_N), z^N(t_N)+\delta_N)) \ge \varepsilon.
\end{equation*}
That is, there is a uniformly positive amount of interactions in an arbitrarily small neighbourhood of $(z(t),t)$, which is not possible (see \cite[Section 4]{BL99}). Hence \eqref{eq:52} holds and therefore, for $N$ large enough,
\begin{equation*}
    |v^N(z^N(t)+,t) - v^+(t)| \le |v^N(z^N(t)+,t) - v^N(\bar{z}^N,t)| + |v^N(\bar{z}^N,t) - v^+(t)| \le 2\varepsilon,
\end{equation*}
which proves \eqref{eq:51}. 

Similarly, we will have, a.e. $t\in [t_0,T]$,
\begin{equation}\label{eq:53}
    v^N(z^N(t)-,t) \to v(z(t)-,t) \quad \text{as } N\to\infty.
\end{equation}

From \eqref{eq:51} and \eqref{eq:53}, there exists a sequence $\delta_N \to 0$ as $N\to \infty$ such that a.e. $t\in [t_0,T]$,
\begin{equation*}
    \dot{z}^N \in [V((z(t)),t)]^{\delta_N}.
\end{equation*}
In other words, $z^N$ is a $\delta_N$-solution of 
\begin{equation} \label{inclus}
    \dot{z}(t) \in V(z(t),t).
\end{equation}
Hence, thanks to \eqref{uniform_conv_2}, we have a uniformly convergent sequence of $\delta_N$-solution $z^N(t)$ of the inclusion \eqref{inclus}. As a result, $z(t)$ is also a solution of this inclusion, thanks to Lemma \ref{lem:Fil-approx}.

Finally, since by Theorem \ref{thm:uniq_trajectory}) the Filippov solution to \eqref{e:LagS} is unique, the whole sequence $z^N$ must converge to $z$. Indeed, if there exists a subsequence $z^{N_l}$ of $z^N$ such that
\begin{equation*}
    z^{N_l} \to y \ne z, \quad \text{ (strongly) uniformly in } [t_0,T],
\end{equation*}
then, by extracting further subsequences if needed, the arguments above show that $y$ also satisfies \eqref{e:LagS} in the sense of Filippov. This contradicts the conclusion of Theorem \ref{thm:uniq_trajectory}.

The proof is complete for the case of front tracking approximations. For the case where $\{v^N\}$ is a sequence of exact solutions, one can first approximate each $v^N$ with a sequence of front tracking approximations $v^{N,n}$, and then use diagonalisation arguments by working with a suitable subsequence $v^{N,n(N)}$.
\end{proof}

We now move on to investigate the approximate Filippov solutions of \eqref{e:LagS} given by the vanishing viscosity approximations of $v$. Consider the parabolic equation 
\begin{equation}\label{scl_viscous}
    v_t + f(v)_x = \epsilon v_{xx},
\end{equation}
where $\epsilon$ is some small positive number. A motivation to study this approximation, apart from the fact that it provides the entropy solution to the original equation \eqref{scl-2} at the zero limit of $\epsilon$, is that it represents a model that takes the deceleration rate of the vehicle before a shock into account. Consider the LWR traffic flow model (see \eqref{eq:traffic}) for example. On the highway, one expects that instead of changing the speed abruptly, 
the driver would slow down when they see increased (relative) density of cars ahead. 
The velocity hence can be written as
    $\tilde{w}(\rho) = w(\rho) - \epsilon\frac{\rho_x}{\rho}$.
The new traffic flow model then is as follows
\begin{equation}\label{traffic_viscous}
    \rho_t + [\rho \tilde{w}(\rho)]_x = \epsilon \rho_{xx}.
\end{equation}
Nevertheless, in general, one may consider the viscous scalar conservation law \eqref{scl_viscous} with an artificial diffusive term on the right-hand side. This artificial diffusive term may be devoid of any physical reasoning but just for the sake of analytical or computational convenience. 

Since \eqref{scl_viscous} is a parabolic equation, the Cauchy problem \eqref{scl_viscous} coupled with some bounded initial value $v_0$ always provides a unique smooth solution $v^\epsilon$. Now consider the trajectory $z^\epsilon$ of a particle starting from $x_0$ and moving along the flow. Assume that the speed of the flow $w$ is a smooth function of $v^\epsilon$. Then it follows that
\begin{equation} \label{trajectory_vis}
    \dot{z}^\epsilon(t) = w(v^\epsilon(z^\epsilon(t),t)), \quad z^\epsilon(t_0) = x_0.
\end{equation}
We investigate the behaviour of $z^\epsilon$ when $\epsilon$ is small and compare it with the trajectory $z$ of the traffic flow, given by (in the sense of Filippov)
\begin{equation} \label{eq:traject}
    \dot{z}(t) = w(v(z(t),t)), \quad z(t_0) = x_0,
\end{equation}
which has been studied earlier. Our aim is to establish a convergence result of $z^\epsilon$ to $z$ as $\epsilon$ goes to 0. We make the following assumption on the speed of the shocks of the original system. 

\begin{assumption}\label{Ass:shockspeed}
The shock speed is always smaller or greater than the speeds of the left and right flows at the shock, that is, either
\begin{equation*}
    s < \min\{w(v_l), w(v_r)\} \quad \text{ or  } \quad s > \max\{w(v_l), w(v_r)\},
\end{equation*}
where the shock speed $s$ is defined as \eqref{Rankine-Hugoniot}, with $v_l$ and $v_r$ are the left and right limits of $v$ at the shock and the flow velocity $w$ satisfy Assumption \ref{as:genw}.
\end{assumption}

The following theorem provides a convergence result for $z^\epsilon$.

\begin{theorem}\label{thm:cont_vis}
Assume that Assumption \ref{Ass:shockspeed} holds at every shock curve of the entropy solution $v$ to \eqref{scl_2} with initial data $v_0 \in L^1(\mathbb{R}) \cap L^\infty(\mathbb{R})$. For each $\epsilon>0$, let $z^\epsilon$ be the solution of \eqref{trajectory_vis} where $v^\epsilon$ is the solution of the viscous scalar conservation law \eqref{scl_viscous} with $v^\epsilon(\cdot,0) = v_0$. Then $z^\epsilon(\cdot)$ converges strongly almost everywhere as $\epsilon \to 0$ to the Filippov solution $z(\cdot)$ of \eqref{eq:traject}.
\end{theorem}

\begin{proof}
We first observe that, since $|\dot{z}^\epsilon| \le 1$, the Arzela-Ascoli theorem will ensure that there exists some absolutely continuous function $z: [0,\infty) \to \mathbb{R}$ such that, up to a subsequence,
\begin{equation} \label{eq:5.9}
    z^\epsilon(\cdot) \to z(\cdot), \quad \text{uniformly in } L^\infty[t_0,T],
\end{equation}
for every $T>t_0$. 

We prove that $z$ is a Filippov solution to \eqref{eq:traject}. 
Thanks to Lemma 6.3.3 in \cite{Daf16}, $\{v^\epsilon(\cdot,t)\}$ is equicontinuous in average. Therefore by Kolmogorov-Riesz theorem (see \cite[Theorem 4.26]{Brz10} or \cite{HH10}), $\{v^\epsilon(\cdot,t)\}$ lies in a compact set of $L^1_{\rm loc}(\mathbb{R})$. Hence, we can find a subsequence, still denoted by $\{v^\epsilon(\cdot,t)\}$, that converges to $v(\cdot,t)$ uniformly. Passing if necessary to a further subsequence, we have that
\begin{equation}\label{eq:5.10}
    v^\epsilon(x,t) \to v(x,t) \quad \text{boundedly almost everywhere on } \mathbb{R} \times [t_0,\infty). 
\end{equation}
Now thanks to Assumption \ref{Ass:shockspeed}, for almost every $t \in [t_0,\infty)$, the point $(z(t),t)$ is a continuity point of $v$. This, together with \eqref{eq:5.10}, ensures that 
\begin{equation}
\label{eq:5.13}
    v^\epsilon(z(t),t) \to v(z(t),t) \quad \text{boundedly almost everywhere on } [t_0,\infty).
\end{equation}
Again, thanks to Lemma 6.3.3 in \cite{Daf16}, 
\begin{equation*}
\int_\mathbb{R} |v^\epsilon(x + y^\epsilon,t) - v^\epsilon(x,t)| dx \le w(|y^\epsilon|) \to 0,
\end{equation*}
as $\epsilon \to 0$, uniformly for $t \in [t_0,T]$. Then up to a subsequence
\begin{equation*}
|v^\epsilon(x + y^\epsilon,t) - v^\epsilon(x,t)| \to 0,
\end{equation*}
as $\epsilon \to 0$, almost everywhere on $\mathbb{R}\times [t_0,T]$. Letting $x = z(t)$ where $(z(t),t)$ is a point of continuity of $v$ and
$$ y^\epsilon = z^\epsilon(t) - z(t),$$
we have, as $\epsilon \to 0$,  
\begin{equation}\label{eq:5.13b}
	|v^\epsilon(z^\epsilon(t),t)-v^\epsilon(z(t),t)| \to 0.
\end{equation}
Finally by writing
\begin{equation}
    |v^\epsilon(z^\epsilon(t),t) - v(z(t),t)| \le |v^\epsilon(z^\epsilon(t),t)-v^\epsilon(z(t),t)| + |v^\epsilon(z(t),t) - v(z(t),t)|,
\end{equation}
and the estimates \eqref{eq:5.13} and \eqref{eq:5.13b}, we conclude that
\begin{equation}
    v^\epsilon(z^\epsilon(t),t) \to v(z(t),t), \quad \text{ almost everywhere in } [t_0,T],
\end{equation}
for every $T>t_0$. This ensures that, for every $\epsilon > 0$ there exists $\delta_{\epsilon}>0$, such that $\delta_{\epsilon} \to 0 $ as $ \epsilon \to 0$, and
\begin{equation}
    \dot{z}^\epsilon(t) \in [v(z(t),t)-\delta_{\epsilon}, v(z(t),t)+\delta_{\epsilon}],
\end{equation}
almost everywhere. Thanks to \eqref{eq:5.10} and Lemma \ref{lem:Fil-approx}, $z$ is a Filippov solution to \eqref{eq:traject}. Due to the uniqueness of $z$, the convergence \eqref{eq:5.9} applies to the whole sequence. This completes the proof.
\end{proof}

\section{H\"older stability of vehicle trajectories of the traffic flow}
\label{sec:stab_Lagrangian}

We shall focus our attention in this section on the traffic flow. \label{subsec:traffic} 
This is one of the most popular applications of the scalar conservation law in one space dimension. We switch to the traditional notation and use $\rho (x,t)$ to denote the quantity of interest in this case, which is the vehicle density (car density) at some given space and time. The velocity of the vehicle flow is still denoted by $w$. The LWR model for traffic flow is derived under general assumptions that the vehicle length is negligible, the road is flat and has only one lane and overtaking is not allowed (see \cite{LW55} and \cite{Ric56}). The Cauchy problem for the LWR model is as follows
\begin{equation} \label{traffic_flow}
    \partial_t\rho + \partial_x (\rho w(\rho)) = 0,
\end{equation}
\begin{equation} \label{traffic_flow_ini}
    \rho(x,0) = \rho_0 (x).
\end{equation}
Denote by $\rho_{max}$ and $w_{max}$ the maximum density and maximum speed of the traffic. We make the following assumption on the car speed and the flux.
\begin{assumption}\label{Ass:carspeed}
The car speed $w: [0, \rho_{max}] \to [0, w_{max}]$ is Lipschitz continuous 
and strictly decreasing function of $\rho$ with $w(\rho_{max}) = 0$. The flux $f(\rho) = \rho w(\rho)$ is a strictly concave function.
\end{assumption}
This assumption is reasonable as we expect the car to go at its maximum speed when there are only a few cars on the road, and to slow down when the car density increases. A typical example for $w$ is that $w$ depends linearly on $\rho$,
\begin{equation}\label{traffic_velo}
    w(\rho) = w_{max}\left( 1 - \frac{\rho}{\rho_{max}} \right).
\end{equation}
By scaling we can assume that $w_{max}=1, \rho_{max}=1$. 
Note that if we replace $\rho$ by $1 - w(\rho)$ in \eqref{traffic_flow}-\eqref{traffic_flow_ini} and denote $v(x,t):= w(\rho(x,t))$ then $v$ follows the following conservation law
\begin{equation} \label{traffic_flow_v}
    v_t + [v(v-1)]_x = 0,
\end{equation}
\begin{equation}\label{traffic_flow_ini_v}
    v(x,t_0) = v_0(x),
\end{equation}
with now a strictly convex flux $\tilde{f}(v)=v(v-1)$ and initial data $v_0 = 1 - \rho_0$. We can see that working with \eqref{traffic_flow_v}-\eqref{traffic_flow_ini_v} is not less general than working with \eqref{traffic_flow}-\eqref{traffic_flow_ini} when the car speed is given as \eqref{traffic_velo}. Note also that, by setting
\begin{equation} \label{transform_Burgers_traffic}
    \rho = 1 - v = \frac{1 - \tilde{v}}{2},
\end{equation}
we obtain the familiar Burgers equation 
\begin{equation} \label{Burgers}
    \tilde{v}_t + \left(\frac{\tilde{v}^2}{2} \right)_x = 0, \quad \tilde{v}(x,t_0) = \tilde{v}_0 (x) := 1-2\rho_0(x),
\end{equation}
for $\tilde{w}$, which is not the velocity here but a function of it as described in \eqref{transform_Burgers_traffic}.

For traffic flow where trapping in the queue is not allowed, then shock speed satisfies Assumption \ref{Ass:shockspeed}, as shown by the following lemma. It follows easily from the Rankine–Hugoniot jump condition, yet will be useful in establishing the stability estimates. 
\begin{lemma}\label{Lem:shockspeed}
Assume that the car density before each time a shock happens is always positive, then shocks travel more slowly than the flows right before and after the shocks. Consequently, trajectories do not lie on shock curves except at countably many points.
\end{lemma}
\begin{proof}
From the Rankine–Hugoniot condition, the speed of a shock curve at the point with left limit $\rho_l$ and right limit $\rho_r$ is given as
\begin{equation*}
    \begin{aligned}
        s & = \frac{\rho_l w(\rho_l) - \rho_r w(\rho_r)}{\rho_l - \rho_r} \\
        & = \rho_l\,\frac{w(\rho_l) - w(\rho_r)}{\rho_l - \rho_r} + w(\rho_r) \\
        & = \rho_r\,\frac{w(\rho_l) - w(\rho_r)}{\rho_l - \rho_r} + w(\rho_l) \\
        & < \min\{w(\rho_l), w(\rho_r)\},
    \end{aligned}
\end{equation*}
since $\rho_l, \rho_r > 0$, and $w$ is a strictly decreasing function of $\rho$. 
\end{proof}

\subsection{Stability with respect to changes in initial field}
We now consider the trajectory $z$ of a car passing through some point $x_0$ at $t=t_0>0$ and travelling at speed $w$. From the previous section, $z$ is the unique Filippov solution to
\begin{equation}
\label{lagrangian_5}
    \dot z(t) = w(\rho(z(t),t)), \quad z(t_0) = x_0.
\end{equation}
The aim is to obtain a suitable convergence rate for some approximation of $z$ with respect to changes in upstream density $\rho_0$, for fixed flux $f$. 
In the following, the Lipschitz constant of $w$ is denoted by 
\begin{align*}
L_w:= \sup_{u\neq v} \left| \frac{w(u) - w(v)}{u-v}\right|.
\end{align*}

\begin{theorem} \label{thm:rate_traffic}
Let $0<m_\rho<1$ and $T>0$ be given. Suppose that $\rho$ and $\bar\rho$ are solutions of  \eqref{traffic_flow}-\eqref{traffic_flow_ini} with initial data $\rho_0, \bar{\rho}_0\in L^1 \cap BV (\mathbb{R};[m_\rho,1])$ respectively, satisfying
\begin{equation}\label{eq:ini_err_rho}
    \|\rho_0-\bar{\rho}_0\|_{L^1\cap L^\infty([-2L_wT,3L_wT])} \le \varepsilon.
\end{equation}
Let $z$ and $\bar{z}$ be the corresponding particle trajectories with the same initial position $x_0$ (so they solve \eqref{lagrangian_5} in the sense Filippov). Then 
\begin{equation*}
    \|z-\bar{z}\|_{L^\infty([t_0,T])} \le  C_\rho\,{\varepsilon}^{1/2}, 
\end{equation*}
 with $C_\rho=1+(T-t_0)(1+\frac{2}{m_\rho})L_w+\frac{1}{m_\rho}(\|\rho_0\|_{BV}+\|\bar\rho_0\|_{BV})$.
\end{theorem}

This result, although not a surprise, seems to be the first stability result with an explicit rate of convergence for particle trajectories given in the context of scalar conservation laws, even for traffic flow. Note that no $L^\infty$ stabilities are expected for entropy solutions, except for classical solutions. 

For an $\varepsilon$-perturbation of $\rho_0$, the idea of the proof is to first estimate the error incurred in $z$ by passing through one single shockwave of $\rho$ of size at least $\varepsilon$, and then to show that the sum of such errors, incurred by passing through the shocks in a time interval of finite length $T$, remains of order $\sqrt{\varepsilon}$ (essentially due to the initial density being of bounded variation). This is done for the front tracking approximation of the solution and then generalised using Theorem \ref{thm:trajectory_cont}.

\begin{proof}[Proof of Theorem \ref{thm:rate_traffic}] 

 For a given $N\in\mathbb{N}$, we start with constructing two simple functions $\rho_0^N$ and $\bar\rho_0^N$ approximating $\rho_0$ and $\bar\rho_0$ as follows. We let
\begin{align*}
&A_{N,j}:=\big\{x\in\R:\rho_0(x)\in[\frac{j-1}{2^N},\frac{j}{2^N})\big\}\quad\mbox{and}\quad\\
&\bar A_{N,j}:=\big\{x\in\R:\bar\rho_0(x)\in[\frac{j-1}{2^N},\frac{j}{2^N})\big\}
\end{align*}
and define
\begin{align*}
\rho_0^N=\sum_{j=1}^{2^N}\frac{j-1}{2^N}\mathbbm{1}_{A_{N,j}}\quad\mbox{and}\quad
\bar\rho_0^N=\sum_{j=1}^{2^N}\frac{j-1}{2^N}\mathbbm{1}_{\bar A_{N,j}}.
\end{align*}
We choose $N>\hat N$ with $\hat N$ large enough so that $1/2^N<\gep<1$ and $\|\rho_0^N-\bar\rho_0^N\|_{L^1\cap L^\infty(-2L_wT,3L_wT)}\le 2\gep$.
 We use the front tracking method to construct approximations of $\rho$, denoted by $\rho^N$ and $\bar\rho^N$ corresponding to $\rho_0^N$ and $\bar\rho_0^N$. That is, we consider $f(\rho)=\rho w(\rho)$ to be approximated by a piecewise linear function $f^N$ whose graph is inscribed by graph of $f$ and $f(j/2^N)=f^N(j/2^N)$ for all $j\in\{0,\dots,2^N\}$. Then, $\rho^N$ and $\bar\rho^N$ are solutions of
\begin{align} \label{e:zeta}
 \zeta_t+(f^N(\zeta))_x=0
\end{align} 
with $\zeta(x,0)=\rho_0^N(x)$ and $\zeta(x,0)=\bar\rho_0^N(x)$ respectively.
 We then define $z^N$ and $\bar z^N$ as 
$$
\dot z^N=w(\rho^N(z^N,t))\quad\mbox{and}\quad\dot {\bar z}^N=w(\bar\rho^n(\bar z^N,t)),
$$
with $z^N(t_0)=\bar z^N(t_0)=x_0$. We first find an upper bound for $\|z^N-\bar z^N\|_{L^\infty(t_0,T)}$.\vskip.2cm

\noindent{\sc Step 1} (Convergence rate for front tracking approximations). In this part, for notational convenience, we drop the superscript $N$ in $\rho^N$, $\bar\rho^N$, $z^N$ and $\bar z^N$. Consider $0<\beta<1$ and
 $t_1>t_0$ to be the smallest time at which
 $$
 \lim_{t\to t_1^+}|\rho(z(t),t)-\bar\rho(\bar z(t),t)|>\gep^\beta.
 $$
 This implies that at $t_1$, one of the particles, which without loss of generality we assume to be $z$, coincides with a shock point. We call this $a_1$ and denote the closest next shock point at this instance to $\bar z$ by $\bar a_1$. We then let the time instance where $\bar z$ hits $\bar a_1$ to be $\tau_1>t_1$. 
 Define
 \begin{align*}
& \iota_1= \lim_{t\to t_1^+} \rho(z(t),t)-\lim_{t\to t_1^-}\rho(z(t),t),\quad\mbox{and}      \\
 &\bar\iota_1= \lim_{t\to t_1^+} \bar\rho(\bar z(t),t)-\lim_{t\to t_1^-}\bar\rho(\bar z(t),t).
 \end{align*}

i) Let $\iota_1>0$ (that is when we have an up-jump at the shock). We consider the Riemann problem around $a_1$ and $\bar a_1$ over the time interval $(t_1,t)$ with $t>\tau_1$ and less than the instant when the next shock is hit by one of the particles.  By Lemma \ref{l:pRiem},  noting that $a_1=z(t_1)$ and $\bar w=w$, we have 
\begin{align}
|\bar z(t)-z(t)|
&\le |w(\bar\rho(\bar z(t),t))-w(\rho(z(t),t)|(t-t_1)\nonumber\\
&\quad+\frac{\iota_1}{m_\rho}|\bar a_1-a_1|+\frac{\lim_{t\to t_1^-}\rho(z(t),t)}{\lim_{t\to t_1^+}\rho(z(t),t)}|\bar z(t_1)-z(t_1)|.\label{e:upjump0}
\end{align}
We then note that $\lim_{t\to t_1^-}\rho(z(t),t)<\lim_{t\to t_1^+}\rho(z(t),t)$ (in the up-jump case), and since
 $$
 \sup_{t_0<t<t_1}|\bar\rho(\bar z(t),t)-\rho(z(t),t)|\le \gep^\beta
 $$
 we have
 $$
 |\bar z(t_1)-z(t_1)|\le (t_1-t_0)L_w\gep^\beta.
 $$
 Furthermore by \cite[Theorem 6.2.3]{Daf16}, 
 $$\|\rho(\cdot,t_1)-\bar\rho(\cdot,t_1)\|_{L^1(0,L_w t_1)}\le \|\rho_0^N-\bar\rho^N_0\|_{L^1(-2L_wT,3L_w T)} \le 2\gep.
 $$ 
 Indeed in the case considered here where $\rho$ and $\bar\rho$ are bounded by 1, constant $s$ in equation (6.2.22) of \cite{Daf16} is bounded by $2L_w$.
 Hence 
 $$
 |a_1-\bar a_1|\gep^\beta\le \|\rho(\cdot,t_1)-\bar\rho(\cdot,t_1)\|_{L^1(0,L_wt_1)}\le 2\gep
 $$
 implying that 
 $$ |a_1-\bar a_1|\le\gep^{1-\beta}.
 $$
We therefore conclude that
\begin{align}
|\bar z(t)-z(t)|
&\le  (t-t_1)L_w \gep^\beta+\frac{\iota_1}{m_\rho}\gep^{1-\beta}+(t_1-t_0)L_w\gep^\beta\nonumber\\
&=  (t-t_0)L_w \gep^\beta+\frac{\iota_1}{m_\rho}\gep^{1-\beta}.\label{e:upjumpi}
\end{align}
The above estimate is valid for $t<t_2$ where $t_2$ is the next time instant at which 
$$\lim_{t\to t_2^+}|\rho(z(t),t)-\bar\rho(\bar z(t),t)|>\gep^\beta.
$$
We also note for the first term in the right-hand side of \eqref{e:upjumpi} we have assumed that over $(\tau_1,t)$,  $\rho$ and $\bar\rho$ are within distance $\gep^\beta$. If this is not the case $\bar a$ can be replaced at the start of the argument by the next shock after which $\rho$ and $\bar\rho$ are within distance $\gep^\beta$. Then a similar argument considering $\lim_{t\to t_1^-}\bar\rho(\bar z(t),t)$ as the value for $\bar\rho$ over $(t_1,\tau_1)$ gives the same estimate.

To derive \eqref{e:upjumpi} we supposed that $z$ hits the shock first, that is $a_1=z(t_1)$. If we have instead $\bar z(t_1)=\bar a_1$, then the only change in \eqref{e:upjumpi} would be the replacement of $\iota_1$ by $\bar\iota_1$.
Hence, in the situation of an up-jump, regardless of which particle hits the shock at $t_1$, we obtain for $\tau_1<t<t_2$
\begin{align}\label{e:upjump1}
 |\bar z(t)-z(t)|\le (t-t_0)L_w \gep^\beta +\frac{\max\{\iota_1,  \bar\iota_1\}}{m_\rho}\gep^{1-\beta}.
 \end{align}
  
 ii) If $\iota<0$, that is when we have a down-jump, it has to be of size $2^{-N}<\gep$. Indeed, for a front tracking solution, all down-jumps after the initial time turn into fans of small shocks, and we note that since $\|\rho_0-\bar\rho_0\|_{L^\infty}<2\gep$, for small enough $\gep$, $t_1>0$. Hence by Lemma \ref{l:pRiem} we have for $\tau_1< t< t_2$
 \begin{align}
 |\bar z(t)-z(t)|&\le |w(\bar\rho(\bar z(t),t))-w(\rho(z(t),t))|(t-t_1)\nonumber\\
&\qquad +\frac{\gep}{m_\rho}\max\big\{a_1-z(t_1), \bar a_1-\bar z(t_1)\big\}+|\bar z(t_1)-z(t_1)|\nonumber\\
&\le (t-t_0)L_w \gep^\beta +\frac{\gep}{m_\rho}\max\big\{a_1-z(t_1), \bar a_1-\bar z(t_1)\big\},\label{e:djump1}
 \end{align}
where we have used the same argument as part (i) to bound $|\bar z(t_1)-z(t_1)|$.

As discuseed above, estimates (\ref{e:upjump1}) and (\ref{e:djump1}) remain valid upto $t_2$. Then at $t_2$ at least one of the particles is at a shock point and the other one is about to hit one, we denote the position of these shocks  at $t_2$ for $z$ and $\bar z$, by $a_2$ and $\bar a_2$ respectively. Defining $\tau_2$, $\iota_2$ and $\bar\iota_2$ similar to above we obtain
\begin{align*}
 |\bar z(t)-z(t)|\le \,&(t-t_0)L_w \gep^\beta+\frac{\gep^{1-\beta}}{m_\rho}\sum_{j=1}^2\max\{\iota_j,  \bar\iota_j,0\}\\
& +\frac{\gep}{m_\rho}\sum_{j=1}^2\max\big\{a_j-z(t_j), \bar a_j-\bar z(t_j)\big\}
\end{align*}
for $t\le t_3$, with $t_3$ the next instant after which $\rho$ and $\bar\rho$ are again more that $\gep^\beta$ apart. Continuing in this fashion until final time $T$ and noting that for a given $\epsilon$, there are finite $k$ number of such shocks ($k$ can depend on $\epsilon$) \cite{BL99}, in the interval $(0,L)$ (with $L<T$), we obtain
\begin{align*}
 |\bar z(T)-z(T)|\le \,&(T-t_0)L_w \gep^\beta+\frac{\gep^{1-\beta}}{m_\rho}\sum_{j=1}^k\max\{\iota_j,  \bar\iota_j,0\}\\
& +\frac{\gep}{m_\rho}\sum_{j=1}^k\max\big\{a_j-z(t_j), \bar a_j-\bar z(t_j)\big\}.
\end{align*}
 By \cite[Theorem 6.2.6]{Daf16} and noting that by Lemma \ref{Lem:shockspeed} $z$ and $\bar z$ do not encounter any of the initial shocks more than once, we have
 $$
 \sum_{j=1}^k\max\{\iota_j,  \bar\iota_j,0\}\le  \sum_{j=1}^k\max\{\iota_j+  \bar\iota_j,0\}\le  \|\rho_0\|_{BV}+\|\bar\rho_0\|_{BV}.
 $$
Moreover, by definition of $t_j$, $a_j$ and $\bar a_j$
 $$
 \sum_{j=1}^k\max\big\{a_j-z(t_j), \bar a_j-\bar z(t_j)\big\}\le 2L_w(T-t_0).
 $$
We hence conclude that
\begin{align*}
\|z-\bar z\|_{L^\infty(t_0,T)}&= |\bar z(T)-z(T)|\\
&\le (T-t_0)L_w \gep^\beta+\frac{\gep^{1-\beta}}{m_\rho}( \|\rho_0\|_{BV}+\|\bar\rho_0\|_{BV})+\frac{2L_w(T-t_0)\gep}{m_\rho}
\end{align*}
as $z$ and $\bar z$ are increasing functions of $t$. The best rate is then achieved when $\beta=1-\beta$, that is $\beta=1/2$.
We hence have, after reinstating superscript $N$,
$$
\|z^N-\bar z^N\|_{L^\infty(t_0,T)}\le C_0\sqrt{\gep},
$$
with $C_0=(T-t_0)(1+\frac{2}{m_\rho})L_w+\frac{1}{m_\rho}(\|\rho_0\|_{BV}+\|\bar\rho_0\|_{BV})$.
\vskip.2cm

\noindent{\sc Step 2} (The general case). Now for any $\rho_0, \bar\rho_0 \in L^1 \cap BV (\mathbb{R};[m_\rho,1])$ such that \begin{equation*}
    \|\rho_0-\bar{\rho}_0\|_{L^1([-2L_wT,3L_wT])} \le \varepsilon,
\end{equation*}
let $\bar{\rho}^N$ and $\rho^N$ be the  the front tracking approximations of $\rho$ and $\bar\rho$, with $\bar{z}^N$ and $z^N$ the corresponding trajectories all constructed as in {\sc Step 1}. We have just proved that,
\begin{equation}\label{eq:4.2}
    \|\bar{z}^N(t) - z^N\|_{L^\infty(t_0,T)} \le C_0\sqrt{\varepsilon}.
\end{equation}
Thanks to Theorem \ref{thm:trajectory_cont}, as $ N\to \infty$ we have
\begin{equation*} \label{eq:4.3}
    \|\bar{z}^N - \bar{z}\|_{L^\infty([t_0,T])} \to 0,
\end{equation*}
\begin{equation*}
    \|z^N - z\|_{L^\infty([t_0,T])} \to 0,
\end{equation*}
for any $t_0>0$. 
Hence, there exist constants $N_1, N_2$ big enough such that
\begin{equation*}
    \|\bar{z}^N - \bar{z}\|_{L^\infty([t_0,T])} \le \frac{\sqrt{\varepsilon}}{2}, \quad\text{for all } N>N_1,
\end{equation*}
\begin{equation*}
    \|z^N - z\|_{L^\infty([t_0,T])} \le \frac{\sqrt{\varepsilon}}{2}, \quad\text{for all } N>N_2.
\end{equation*}
Choose $N > \max\{\hat N, N_1, N_2\}$, together with \eqref{eq:4.2} we obtain
\begin{align*}
   \|\bar{z} - z\|_{L^\infty([t_0,T])} & \le \|\bar{z}^N-\bar{z}\|_{L^\infty([t_0,T])} + \|z^N-z\|_{L^\infty([t_0,T])} + \|\bar{z}^N - z^N\|_{L^\infty([t_0,T])} \\
   & \le (1+C_0)\sqrt{\varepsilon}.
\end{align*}
This completes the proof. 
\end{proof}

\begin{lemma}\label{l:pRiem}
Let $w,\bar w\in\Lambda$, and $f,\bar f:(0,1]\to\R$ be continuous and concave functions.
Let $\zeta$ and $\bar\zeta$ satisfy 
\begin{align*}
 &\zeta_t + (f(\zeta))_x = 0,\hspace{2.65cm}\bar\zeta_t + (\bar f(\bar\zeta))_x = 0,\\
 & \zeta(x,0) = \begin{cases}
\rho_l & \text{ if } x<a \\ 
\rho_r & \text{ if } x>a, 
\end{cases}  
\hspace{1.2cm}
\bar{\zeta} (x,0) = \begin{cases}
\bar \rho_l  & \text{ if } x < \bar a\\ 
\bar \rho_r  & \text{ if } x > \bar a,
\end{cases} 
\end{align*}
respectively and suppose that $f(\rho_i)=\rho_i w(\rho_i)$ and $\bar f(\bar\rho_i)=\bar\rho_i \bar w(\bar\rho_i)$ for $i=l,r$.\\
Let $z:[t_0,T]\to\R$ and $\bar z:[t_0,T]\to\R$ be solutions of
 $$
 \dot z= w(\zeta(z,t))\quad\mbox{and}\quad\dot {\bar z}=\bar w(\bar\zeta(\bar z,t)),
 $$
 respectively with $z(t_0)=z_0<a$ and $\bar z(t_0)=\bar z_0<\bar a$.\\ 
 Then, for any $t>t_0$, 
$$
\bar z(t)-z(t)
=(\bar w(\bar\rho_r)-w(\rho))(t-t_0)+\frac{\bar\rho_r-\bar\rho_l}{\bar\rho_r}(\bar a-\bar z_0)-\frac{\rho_r-\rho_l}{\rho_r}(a-z_0)+\bar z_0-z_0.
$$
\end{lemma}

\begin{proof}
We have
\[ \zeta(x,t) = \begin{cases}
\rho_l & \text{ if } x< a+\lambda t \\ 
\rho_r & \text{ if } x> a+\lambda t. 
\end{cases} \]
with
\begin{equation*}
    \lambda = \frac{f(\rho_l)-f(\rho_r)}{\rho_l - \rho_r} = \frac{\rho_l w(\rho_l)- \rho_r w(\rho_r)}{\rho_l - \rho_r},
\end{equation*}
and
\[ \bar{\zeta} (x,t) = \begin{cases}
\bar \rho_l & \text{ if } x< \bar a+ \bar{\lambda}t \\ 
\bar \rho_r & \text{ if } x> \bar a + \bar{\lambda}t, 
\end{cases} \]
with
\begin{equation*}
    \bar{\lambda} = \frac{\bar\rho_l \bar w(\bar\rho_l)- \bar\rho_r \bar w(\bar\rho_r)}{\bar\rho_l - \bar\rho_r}.
\end{equation*}
The particle $z$ first travels with speed $w(\rho_l)$ until it hits the shock,
then travels with speed $w(\rho_r)$. The hitting time $\tau$ can be calculated as $ w(\rho_l) \tau + z_0 = \lambda \tau+a$, thus
\begin{equation*}
    \tau = \frac{a-z_0}{w(\rho_l)-\lambda}.
\end{equation*}
 The trajectory of $z$ is given by 
\[ z(t) = \begin{cases}
w(\rho_l) (t-t_0)+ z_0 & \text{ if } t-t_0< \tau \\ 
w(\rho_r)(t-t_0- \tau) + w(\rho_l)\tau + z_0 & \text{ if } t-t_0\ge \tau. 
\end{cases} \]
On the other hand, the particle $\bar{z}$ first travels with speed $\bar w(\bar{\rho}_l)$ until it hits the shock and travels with speed $\bar w(\bar{\rho}_r)$ after that. Similarly, 
\begin{equation}\label{eq:5.31}
    \bar z(t) =\begin{cases}
\bar w(\bar{\rho}_l) (t-t_0)+ \bar z_0 & \text{ if } t-t_0< \bar\tau \\ 
\bar w(\bar{\rho}_r)(t -t_0- \bar\tau) + \bar w(\bar{\rho}_l)\bar\tau + \bar z_0 & \text{ if } t-t_0\ge \bar\tau,
\end{cases}
\end{equation}
where $\bar\tau$ is the hitting time
\begin{equation}\label{eq:hitting_2}
    \bar\tau = \frac{\bar a- \bar z_0}{\bar w(\bar{\rho}_l) - \bar{\lambda}}.
\end{equation}
Hence after the particle has passed both shock points, that is $t>\max\{\tau,\bar\tau\}+t_0$, we have, setting $w_i:=w(\rho_i)$ and $\bar w_i:=\bar w(\bar\rho_i)$ for $i=r,l$,
\begin{align*}
\bar z(t)-z(t)&= (\bar w_r-w_r)(t-t_0)+(\bar w_l-\bar w_r)\bar\tau-(w_l-w_r)\tau+\bar z_0-z_0.
\end{align*}
Noting that 
\begin{align*}
(w_l-w_r)\tau=(w_l-w_r)\frac{(\rho_r-\rho_l)(a-z_0)}{\rho_r(w_l-w_r)}=\frac{\rho_r-\rho_l}{\rho_r}(a-z_0)
\end{align*}
and similarly
$$
(\bar w_l-\bar w_r)\bar\tau=\frac{\bar\rho_r-\bar\rho_l}{\bar\rho_r}(\bar a-\bar z_0),
$$
we have
\begin{align*}
\bar z(t)-z(t)
=(\bar w_r-w_r)(t-t_0)+\frac{\bar\rho_r-\bar\rho_l}{\bar\rho_r}(\bar a-\bar z_0)-\frac{\rho_r-\rho_l}{\rho_r}(a-z_0)+\bar z_0-z_0
\end{align*}
and the result follows.
\end{proof}

\subsection{Stability with respect to changes in flux function}
In this section, we study the stability of the particle trajectories with respect to small changes in the flux function (with the initial field being fixed). We define
\begin{equation}\label{Lip_space}
    \Lambda := \left\{w:[0,1] \to \mathbb{R}, w \mbox{ is strictly decreasing}, w(1)=0,\mbox{ and } \|w\|_{\rm Lip} < \infty \right\},
\end{equation} 
where $\|\cdot\|_{\rm Lip}$ is defined as \eqref{Lip-norm}. Let $\bar w\in \Lambda$ and $\bar\rho$ satisfy
\begin{equation} \label{traffic_flowb}
    \partial_t\bar\rho + \partial_x (\bar\rho\, \bar w(\bar\rho)) = 0,
\end{equation}
\begin{equation} \label{traffic_flow_inib}
    \bar\rho(x,0) = \rho_0 (x).
\end{equation}
We have the following stability estimate.

\begin{theorem}\label{thm:rate_flux}
Assume that $\rho$ and $\bar\rho$ are the entropy solutions to (\ref{traffic_flow})-(\ref{traffic_flow_ini}) and (\ref{traffic_flowb})-(\ref{traffic_flow_inib}) respectively with the same initial data $\rho_0 \in L^1 \cap BV (\mathbb{R};[m_\rho,1])$, $0<m_\rho<1$, and 
with $w,\bar w\in\Lambda$ satisfying
\begin{equation*}
    \|w- \bar w\|_{\rm Lip} \le \varepsilon.
\end{equation*}
Then the corresponding particle trajectories $z$ and $\bar{z}$ satisfy
\begin{equation*}
    \|z-\bar{z}\|_{L^\infty([t_0,T])} \le  C_w\sqrt{\varepsilon},
\end{equation*}
for any given $T>0$, and with $C_w=1+2(T-t_0)(1+\frac{2}{m_\rho})\|w\|_{\rm Lip}+\frac{2}{m_\rho}\|\rho_0\|_{BV}$.
\end{theorem}

\begin{proof}
We argue as in the proof of Theorem \ref{thm:rate_traffic} and 
first consider $\rho^N$ and $\bar\rho^N$, the front tracking approximations of $\rho$ and $\bar\rho$,
 and their corresponding trajectories $z^N$ and 
$\bar z^N$ respectively. We follow a similar argument to {\sc Step 1} of proof of  Theorem \ref{thm:rate_traffic} and employ Lemma \ref{l:pRiem}, the only difference here is that in the right-hand side of inequalities (\ref{e:upjump0}) and (\ref{e:djump1})  $w(\bar\rho)$ is replaced with $\bar w(\bar\rho)$ and we write instead
\begin{align*}
|\bar w(\bar\rho(\bar z(t),t))-w(\rho(z(t),t)|&\le | w(\bar\rho(\bar z(t),t))-w(\rho(z(t),t)|+|\bar w(\bar\rho(\bar z(t),t))-w(\bar\rho(z(t),t)|\\
&\le L_w\gep^\beta+\|w-\bar w\|_{L^\infty}\\
&\le  L_w\gep^\beta+\gep,
\end{align*}
as $\rho\le 1$ and hence $ \|w-\bar w\|_{L^\infty}\le \|w-\bar w\|_{\rm Lip}$.
We therefore obtain 
\begin{align*}
|\bar z^N(T)-z^N(T)|\le &(T-t_0)( L_w\gep^\beta+\gep)+\frac{\gep^{1-\beta}}{m_\rho}\sum_{j=1}^k\max\{\iota_j,  \bar\iota_j,0\}\\
& +\frac{\gep}{m_\rho}\sum_{j=1}^k\max\big\{a_j-z(t_j), \bar a_j-\bar z(t_j)\big\}.
\end{align*}
Since, again by  \cite[Theorem 6.2.6]{Daf16} and as by Lemma \ref{Lem:shockspeed} $z$ and $\bar z$ do not encounter any of the initial shocks more than once, we have
$$
\sum_{j=1}^k\max\{\iota_j,  \bar\iota_j,0\}\le \sum_{j=1}^k \max\{\iota_j+  \bar\iota_j,0\}\le 2\|\rho_0\|_{BV},
$$
and noting that
\begin{align*}
\sum_{j=1}^k\max\big\{a_j-z(t_j), \bar a_j-\bar z(t_j)\big\}
&\le \sum_{j=1}^k (a_j-z(t_j)) + (\bar a_j-\bar z(t_j))\\ 
&\le \|w\|_{\rm Lip}(T-t_0)+\|\bar w\|_{\rm Lip} (T-t_0)\\
&\le (2L_w+\gep)(T-t_0).
\end{align*}
we obtain
\begin{align*}
|\bar z^N(T)-z^N(T)|
&\le (T-t_0)( L_w\gep^\beta+\gep)+\frac{2\|\rho_0\|_{BV}}{m_\rho}\gep^{1-\beta}+\frac{\gep}{m_\rho}(2L_w+\gep)(T-t_0)\\
&\le 2(T-t_0)(1+\frac{2}{m_\rho})\gep^{\beta}+\frac{2\|\rho_0\|_{BV}}{m_\rho}\gep^{1-\beta}
\end{align*}
as $\beta<1$, and we get the best rate in $\gep$ when we choose $\beta=1/2$. We hence have
$$
\|\bar z^N-z^N\|_{L^\infty}\le C_0\sqrt{\gep}
$$
with $C_0=2(T-t_0)(1+\frac{2}{m_\rho})+\frac{2\|\rho_0\|_{BV}}{m_\rho}$. Then, the same argument as {\sc Step 2} of Theorem \ref{thm:rate_traffic} gives the result.
\end{proof}


\section{Bayesian inverse problems for initial field and flux function}
\label{sec:Bay_recovery}
\subsection{Bayesian inverse problems}
As we have discussed in the introduction, \label{subsec:bip}
a mathematical formulation of the inverse problem of finding the unknown $u$ from the observed data $y$ reads as
\begin{equation}
\label{eq:IP}
    y = G(u),
\end{equation}
where $G$ is the observation operator. The Bayesian approach for inverse problems typically starts with the observation that the data $y$ is usually perturbed by noise, hence a more appropriate formulation for \eqref{eq:IP} should be 
\begin{equation}
\label{eq:BIP}
y = G(u) + \xi,
\end{equation}
where \(\xi\) denotes the \emph{observable noise} \cite{Stu10, KS05}. It is then natural to treat the data $y$ and the unknown $u$ as random variables. The solution of the inverse problem \eqref{eq:BIP} will also be a random variable, denoted by \(u|y\) (\(u\) given \(y\)). In the Bayesian framework, we store our prior information about $u$ (before the data \(y\) is taken into account) as a probability distribution and use Bayesian inference to calculate the posterior $\mu^y$ which stores, in turn, the information about \(u|y\). In other words, the Bayesian approach provides a probability distribution \(\mu^y\) carrying information about \(u\), instead of finding the exact solution \(u\) which is not possible in most ill-posed problems.



The common framework for Bayesian inverse problems (BIP) is as follows. Assume that the unknown \(u\) lies in \(X\), the data \(y\) is given in \(Y\), with \(X, Y\) Banach spaces. We suppose that a prior probability measure $\mu_0$ on $(X,\mathcal{B}(X))$ is given and $\xi\sim\mathbb{Q}_0$ with $\mathbb{Q}_0$ known. We assume also that $\xi$ and $u$ are independent with $y|u\sim\mathbb{Q}_u$ which for some $\Phi:X\times Y\to \R$ (obtained through equation (\ref{eq:BIP})) satisfies
$$
\frac{{\rm d}\mathbb{Q}_u}{{\rm d}\mathbb{Q}_0}=\exp(-\Phi(u;y)).
$$
Let the measure $\nu$ be defined by $\nu(\ud u,\ud y):=\mu_0(du)\mathbb{Q}_u(dy)$. Then the following version of Bayes' theorem is proved in \cite{DS16}.
\begin{theorem}\label{t:GenBay} Assume that
$\Phi:X\times Y \to \R$ is $\mu_0\otimes\mathbb{Q}_0$ measurable and that, for $y$
$\mathbb{Q}_0$-a.s.,
\begin{equation}
\label{eq:mev}
Z:=\int_{X}\exp\big(-\Phi(u;y)\big)\mu_0(du)\in (0,\infty).
\end{equation}
Then the conditional distribution of $u|y$ exists under $\nu$, and is denoted by
$\mu^y$. Furthermore $\mu^y\ll\mu_0$ and, for $y\in Y$, it holds $\nu$-a.s.,
\begin{equation}\label{eq:Bay}
\frac{d\mu^y}{d\mu_0}(u)=\frac{1}{Z}\exp\big(-\Phi(u;y)\big).
\end{equation}
\end{theorem}

When the data space $Y$ is finite-dimensional and the observational noise $\xi$ is non-degenerate Gaussian, the measure $\mu^y$ is continuous in $y$ under very mild conditions. A proof of the following result can be found in \cite{Lat20}.
\begin{theorem}\label{t:wpBay}
Let $Y:=\R^{J}$, $\xi$ have Gaussian distribution $\mathcal{N}(0,\Gamma)$ with $\Gamma$ symmetric positive definite, and $G:X\to\R$ be $\mu_0$-measurable. Then for any $y\in Y$ and $\{y^k\}\subset Y$ with $y^k\to y$ as $k\to\infty$, we have
$$
d_{{\rm Hell}}(\mu^y,\mu^{y^k})\to 0\quad\mbox{as}~~~k\to\infty
$$
where $d_{\mathrm{Hell}}(\mu,\mu')$ is the \emph{Hellinger distance} between two measures \(\mu\) and \(\mu'\),
\begin{equation*}
\big(d_{\mathrm{Hell}}(\mu,\mu')\big)^2 := \frac{1}{2}\int_X\left(\sqrt{\frac{d\mu}{d\mu_0}} - \sqrt{\frac{d\mu'}{d\mu_0}}\right)^2 d\mu_0.
\end{equation*}
\end{theorem}

In practice, solving a PDE normally involves some sort of approximation. The question of whether the perturbed posterior arising from the approximation of the forward model converges to the posterior given by the exact one is then importantly necessary to study. Let \(\mu\) (we drop the variable \(y\) since it does not have any explicit role in this task) be the solution of the Bayesian inverse problem \eqref{eq:BIP}, which is given by 
\begin{equation} \label{posterior}
\dfrac{d \mu}{d\mu_0} (u) = \dfrac{1}{Z} \exp (-\Phi(u)),
\end{equation}
\begin{equation} \label{normalization}
Z = \int_X \exp (-\Phi(u))\mu_0(du).
\end{equation}
Let \(\mu^N\) be the measure defined by
\begin{equation} \label{posterior-N}
\dfrac{d \mu^N}{d\mu_0} (u) = \dfrac{1}{Z^N} \exp (-\Phi^N(u)),
\end{equation}
\begin{equation} \label{normalization-N}
Z^N = \int_X \exp (-\Phi^N(u))\mu_0(du),
\end{equation}
where \(\Phi^N(u)\) is some approximation of \(\Phi(u)\). The existence of $\mu$ and $\mu^N$ are guaranteed by Theorem \ref{t:GenBay}. If the data is finite and the noise is Gaussian, $\xi \sim \mathcal{N}(0,\Gamma)$, then $\Phi$ and $\Phi^N$ can be defined as
\begin{equation}\label{eq:Phi_PhiN}
    \Phi (u) = \frac{1}{2}|y - G(u)|_{\Gamma}^2, \quad \Phi^N (u) = \frac{1}{2} |y - G^N(u)|_{\Gamma}^2,
\end{equation}
where $G^N(u)$ is some approximation of $G(u)$. We would like to see whether (and how) the approximation \(G^N(u)\) of \(G(u)\) translates to the approximation \(\mu^N\) of \(\mu\). The following theorem which follows from \cite[Theorem 4.9]{DS16}, gives a sufficient condition for that translation.

\begin{theorem}\label{thm:app_gen}
Assume that there exists a measurable function \(M(\cdot): \mathbb{R}_+\to \mathbb{R}_+\) such that \(G(u), G^N(u): X \to \mathbb{R}^J\) satisfy the following condition for all \(u\in X\),
\begin{equation} \label{eq:cont_FWmap}
|G(u) - G^N(u)| \le M(\|u\|_X)\psi(N), 
\end{equation}
where \(\psi(N)\to \infty\) as \(N\to \infty\).
Suppose in addition that $\mu_0$ is a probability measure on $X$ such that
\begin{equation*} 
\exp\left(M(\|u\|_X)\right) \in L^1_{\mu_0}(X;\mathbb{R}).
\end{equation*}
Then, there exists $C>0$ such that
\[d_{\mathrm{Hell}}(\mu,\mu^N) \le C\psi(N),\]
for all $N$ sufficiently large.
\end{theorem}

When the data comes from discrete measurements of the particle trajectories, 
we show the consistency of the posterior with respect to appropriate approximations 
of solution of  \eqref{scl}-\eqref{ini}. 
For the case that observations are pointwise evaluations of the solution $v$ of  \eqref{scl}-\eqref{ini}, 
the development of shockwaves prevents such a continuity result for the approximations.
It is however possible to obtain some partial results. We leave this to Section \ref{sec:app_Euler} for further discussion.

In the rest of this section we consider $y\in\R^J$, $\Gamma$ a diagonal matrix with nonzero members all equal $\gamma^2$, and $\mu$ and $\mu^N$ satisfying
\begin{align}
 &  \frac{d\mu}{d\mu_0}(u)=\frac{1}{Z}\exp\big(-\frac{1}{2\gamma^2}|y-G(u)|^2\big)  \label{e:posterior}\\
 &   \frac{d\mu^N}{d\mu_0}(u)=\frac{1}{Z^N}\exp\big(-\frac{1}{2\gamma^2}|y-G^N(u)|^2\big)\label{e:Nposterior}
\end{align}
with $Z$ and $Z^N$ given as in \eqref{normalization} and \eqref{normalization-N}.

\subsection{BIP for initial field with discrete measurements of particle trajectories}
\label{subsec:app_rate}


We suppose in this section that we make a finite number of noisy observations $y=(y_1~y_2~\dots~y_J)^T$ from one particle trajectory,
\begin{equation}\label{eq:LagData}
y_j=z(t_j)+\xi_j,\quad~~j=1,\dots ,J,    
\end{equation}
with $z$ defined as the Filippov solution of \eqref{e:LagS}. We suppose that $\xi_j\sim \mathcal{N}(0,\gamma^2)$ and are independent. 

Here we are interested in recovering the upstream field. We let
\begin{align}\label{e:L1G(u)}
    \begin{array}{l}
         u := v(\cdot,0) \\
         G(u):=(G_1(u),~\dots~, G_J(u)),\mbox{~ with ~}G_j(u)=z(t_j).
    \end{array}
\end{align}
where $t_J>t_{J-1}>\dots>t_1>0$. We have the following well-posedness result.
\begin{theorem}\label{t:LwpBay}
Suppose $y\in\R^J$ is given by \eqref{eq:LagData}. Let $\mu_0(X)=1$ where $X = L^1(\mathbb{R})\cap BV(\mathbb{R})$. Then  the posterior measure $\mu^y$ given by 
\begin{equation*}
    \frac{d\mu^y}{d\mu_0}(u)=\frac{1}{Z}\exp\big(-\frac{1}{2\gamma^2}|y-G(u)|^2\big),
\end{equation*}
with  $Z=\int_{X} \exp\big(-\frac{1}{2\gamma^2}|y-G(u)|^2\big)\,\ud\mu_0$ is  well-defined and continuous in Hellinger distance with respect to $y$.
\end{theorem}

\begin{proof}
We write $G=\mathcal{O}\circ S_l$ with $\mathcal{O}:X\to\R^J$ the point observation operator $z\mapsto y$ and $S_l$ the forward operator mapping $u\in X$ to $z\in C_b([t_0,T])$ with $C_b$ denoting the space of bounded continuous functions. 
Consider $\{u^N\}\subset L^1\cap L^\infty$ and suppose that $u^N\to u$ in $L^1$ and let $v^N$ be the solution of the conservation with $v^N(0)=u^N$. Then we have $v^N\to v$ in $L^1(0,T; L^1\R)$. 
Let $z^N:[t_0, T]\to\R$, $t_0>0$, be the unique Filippov solution to 
$$
\dot{z}^N(t) \in [w(v^N(z^N(t)+,t)), w(v^N(z^N(t)-,t))],~~~ \mbox{ with } z^N(t_0)=z(t_0)=x_0.
$$
By Theorem \ref{thm:trajectory_cont} we have $z^N\to z$ uniformly. This concludes the continuity of $S_l:u\mapsto z$. Since the point observation operator $\mathcal{O}$ is continuous we have $G=\mathcal{O}\circ S_l:X\to\R^J$ is continuous. It is evident that $Z<1$ and since $G$ is bounded we have $Z>0$. The result follows by Theorem \ref{t:GenBay} and \ref{t:wpBay}. 
\end{proof}


We now investigate the continuity property of the posterior $\mu^y$ given in Theorem \ref{t:LwpBay} with respect to perturbations in the forward problem. Let $\mu$ and $\mu^N$ satisfy \eqref{e:posterior} and \eqref{e:Nposterior} respectively with $G$ given as \eqref{e:L1G(u)} and 
\begin{equation}\label{e:L1G_N(u)}
    G^N(u) = \{z^N(t_j)\}_{j=1}^J,
\end{equation}
where $z^N(t)$ denotes some approximation of the particle trajectory $z(t)$. 

The first convergence result concerns the approximate measure $\mu^N$ arising in approximations of the initial condition of the forward problem (giving a sequence of exact solutions) or front tracking approximations. 

\begin{theorem}\label{thm:p_conv-0-rate}
Let $\{v^N\} \subset L^1\cap BV(\mathbb{R})$ be a sequence of exact solutions (or front tracking approximations) that converges in $L^1$ to the entropy solution $v$ of \eqref{scl}-\eqref{ini}, and consider $\{z^N\}$ and $z$ to be the corresponding trajectories respectively. Let $\mu^N, \mu$ be given as in \eqref{e:posterior}-\eqref{e:Nposterior} where $G, G^N$ are defined by \eqref{e:L1G(u)}-\eqref{e:L1G_N(u)}. Then
\begin{equation}\label{eq:conv_hell}
    d_{\mathrm{Hell}}(\mu^N,\mu) \to 0,
\end{equation}
as $N\to \infty$.
\end{theorem}

\begin{proof}
From the definition of Hellinger distance, the formulations \eqref{posterior}-\eqref{posterior-N} and the basic inequality $(a+b)^2 \le 2(a^2 + b^2)$ we have
\begin{equation*} 
\begin{aligned} 
d_{\rm Hell}(\mu,\mu^N)^2 & = \frac{1}{2}\int_X\left(\sqrt{\frac{d\mu}{d\mu_0}} - \sqrt{\frac{d\mu^N}{d\mu_0}}\right)^2 \mu_0(du)\\ 
& \le I_1 + I_2,
\end{aligned}
\end{equation*}
where 
\begin{equation}\label{eq:I1} 
I_1 = \frac{1}{Z} \int_X \left(e^{-\frac{1}{2}\Phi(u)} - e^{-\frac{1}{2}\Phi^N(u)}\right)^2 \mu_0(du),
\end{equation}
\begin{equation}\label{eq:I2}
I_2 = \left|\frac{1}{\sqrt{Z}}-\frac{1}{\sqrt{Z^N}}\right| \int_X e^{-\Phi^N(u)} \mu_0(du).
\end{equation}
Since the noise is Gaussian, it follows from \eqref{eq:Phi_PhiN} that $$|\Phi(u)-\Phi^N(u)| \le C|G(u) - G^N(u)|,$$ 
where $C$ is a constant independent of $u$. 
Together with the locally Lipschitz property of $e^{-x}$ for $x>0$, we have 
\begin{equation*}
\begin{aligned}
I_1 &\le C\int_X \left|\Phi(u)-\Phi^N(u)\right|^2 \mu_0(du) \\
& \le C\int_X \left|G(u)-G^N(u)\right|^2 \mu_0(du).
\end{aligned}
\end{equation*}
With $G(u) $ and $ G^N(u)$ defined by \eqref{e:L1G(u)} and \eqref{e:L1G_N(u)}, thanks to Theorem \ref{thm:trajectory_cont}, for every $u \in X$, $G^N(u) \to  G(u)$. By dominated convergence theorem, this ensures that $I_1 \to 0$ as $N \to \infty$.

For $I_2$, it is sufficient to show that $|Z - Z^N| \to 0$. We have
\begin{align*} 
|Z-Z^N| &\le \int_X |\exp(-\Phi(u)-\exp(-\Phi^N(u))|\mu_0(du) \\
&\le \int_X |\Phi(u)-\Phi^N(u)|\mu_0(du).
\end{align*}
The same arguments as above lead to $|Z - Z^N| \to 0$, therefore $I_2 \to 0$ as $N \to \infty$. This completes the proof of Theorem \ref{thm:p_conv-0-rate}.
\end{proof}

\begin{remark} \label{rm_wass_con}
    The use of Hellinger distance is for convenience and to be consistent with most of the works in literature. Other choices of distances between measures are possible. We mention here an important one, the Wasserstein distance, whose usage has been rising recently in the theory of optimal transport, statistics and machine learning (see for instance, \cite{PC19}). A version of the convergence \eqref{eq:conv_hell} in Wasserstein distance (of the first order) $d_{\mathrm{Wass}}$ can be proved. Indeed, thanks to the famous Kantorovich–Rubinstein duality theorem (see \cite{Vil08}), we may write 
\begin{equation*} 
\begin{aligned}
d_{\mathrm{Wass}} (\mu, \mu^N) & =  \sup_{\|h\|_{\rm Lip} \le 1}  \left| \int_X h d\mu - \int_X h d\mu^N \right| \\
& = \sup_{\|h\|_{\rm Lip} \le 1} \left|\int_X h(u)\frac{\exp(-\Phi(u))}{Z}  - \frac{\exp(-\Phi^N(u))}{Z^N} \mu_0(du)\right|,
\end{aligned}
\end{equation*}
and then use the estimate
\begin{equation*}
    |h(u)| \le \|h\|_{\rm Lip}\|u\|_X + |h(0)| \le C (1+\|u\|_X) 
\end{equation*}
to get rid of $h$. The rest of the proof can be carried out similarly, with the help of Fernique's theorem where appropriate. See also \cite{Spr20}.
    In the rest of the paper, all convergence results in Hellinger distance also apply to Wasserstein distance, with suitable modifications as we discussed above.
\end{remark}

The convergence result for the approximate posterior when one uses the vanishing viscosity approximation for the forward problem follows next. 
The proof is entirely similar to the one of Theorem \ref{thm:p_conv-0-rate} 
so will be omitted here.

\begin{theorem}\label{thm:p_conv-0-rate_vis}
Assume that Assumption \ref{Ass:shockspeed} holds at every shock curve of the entropy solution $v$ to the equation \eqref{scl_2}. Let $v^\epsilon$ be the vanishing viscosity approximation of $v$ and $z^\epsilon$ be the corresponding solution to \eqref{trajectory_vis}. Let $\mu$ be given as in \eqref{e:posterior} with $G$ given in \eqref{e:L1G(u)}. Assume that $\mu^\epsilon$ is defined as follows
\begin{equation*}
    \frac{d\mu^\epsilon}{d\mu_0}(u)=\frac{1}{Z^\epsilon}\exp\big(-\frac{1}{2\gamma^2}|y-G^\epsilon(u)|^2\big),    
\end{equation*}
where $G^\epsilon(u) = \{z^\epsilon(x_j,t_j)\}_{j=1}^J$.
Then $\mu^\epsilon$ converges to $\mu$ in the sense that 
\begin{equation}
    d_{\rm Hell}(\mu^\epsilon,\mu) \to 0,
\end{equation}
as $\epsilon\to 0$.
\end{theorem}


\subsection{The case of traffic flow}
In this section we study the convergence properties of the approximated posterior \eqref{e:Nposterior} when the unknown $u$ is either the upstream density $\rho_0$ or the flux function (equivalently velocity function $w$) of the traffic flow model \eqref{traffic_flow}-\eqref{traffic_flow_ini}.
\subsubsection{Approximat BIP for the initial field}
 Let $G$ defined as in \eqref{e:L1G(u)}. 
We consider the approximate Bayesian inverse problem \eqref{e:Nposterior} where the approximate observation map is given as
\begin{equation}\label{e:TrafG_N(u)}
 G^N(u) := G(u^N),
\end{equation}
with $u^N\in X$ an approximation of $u\in X$ with $X=L^1 \cap BV (\mathbb{R};(0,1))$. 
The well-posedness of $\mu$ and $\mu^N$ is guaranteed by Theorem \ref{t:LwpBay}.
The following theorem
provides a convergence rate for the approximation of the posterior in terms of the upstream density.

\begin{theorem}\label{thm:conv_traffic}
Let $X=L^1 \cap BV (\mathbb{R};(0,1))$. Suppose that for any $u\in X$, the approximating sequence $\{u^N\} \subset X$ satisfies
\begin{equation*}
    \|u^N-u\|_{L^1\cap L^\infty} \le \psi(N), \quad \psi(N) \to 0 \text{ as } N\to \infty.
\end{equation*}
Assume that $\mu_0(X)=1$, and
\begin{align}\label{e:intC}
\int_X \frac{1+\|u\|_{BV}^2}{m_u^2}\;\mu_0(\ud u)<\infty
\end{align}
where $m_u=\inf_{x\in\mathbb{R}}u(x)$.
Let $\mu$ and $\mu^N$ be given as in \eqref{e:posterior} and \eqref{e:Nposterior} with $G^N$ defined by \eqref{e:TrafG_N(u)}.

Then 
\begin{equation*}
    d_{\mathrm{Hell}}(\mu^N,\mu) \le C\sqrt{\psi(N)},
\end{equation*}
as $N\to\infty$.
\end{theorem}

\begin{proof}
By processing similarly as in the proof of Theorem \ref{thm:p_conv-0-rate}, we arrive at
\begin{equation*}
    d_{\rm Hell}(\mu,\mu^N)^2 \le I_1 + I_2,
\end{equation*}
where $I_1$ and $I_2$ are given by \eqref{eq:I1} and \eqref{eq:I2}. We then have
\begin{equation*}
    I_1 \le C\int_X |G(u) - G(u^N)|^2 \mu_0(\ud u) \le C\psi(N)\int_X \frac{1+\|u\|^2_{BV}}{m_u^2}\mu_0(\ud u)\le C\psi(N) ,
\end{equation*}
by (\ref{e:intC}) and in the second inequality we have used Theorem \ref{thm:rate_traffic}. A similar argument leads to 
\begin{equation*}
    I_2 \le C|Z - Z^N|^2 \le C\psi(N).
\end{equation*}
The result then follows.
\end{proof}

\begin{remark} A prior satisfying the conditions of the above theorem can be constructed as follows. 
Let $\nu_0$ be a Gaussian measure with $\nu_0(W^{1,1})=1$ where $W^{1,1}$ is the space of integrable functions with integrable derivatives on $\R$. We note that $W^{1,1}(\R)\subset L^1\cap BV(\R)$.
Let $v\sim \nu_0$ and 
$$
u(x)=F(v(x)):=\left\{
\begin{array}{ll}
\frac{1}{2}{\rm e}^{v(x)},&\mbox{if }  v(x)\le 0\\
1-\frac{1}{2}{\rm e}^{-v(x)},&\mbox{if }  v(x)>0.
\end{array}
\right.
$$
Consider $\mu_0:=\nu_0\circ F^{-1}$. We have
\begin{align*}
m_u=\min_{x\in\R} u(x)\ge \frac{1}{2}\,{\rm e}^{-\|v\|_{L^\infty}}
\end{align*}
and
 $$\|u\|_{BV}=\int_\R |u'|\,\ud x\le \frac{1}{2}\int_\R |v'|\,\ud x\le \frac{1}{2}\|v\|_{W^{1,1}}.$$ 
 Hence
 $$
 \int_X \frac{1+\|u\|_{BV}^2}{m_u^2}\;\mu_0(\ud u)\le \int_{W^{1,1}}\|v\|_{W^{1,1}}^2 {\rm e}^{2\|v\|_{W^{1,1}}}\,\ud \nu_0(\ud v)<\infty
$$
since $\|v\|_{L^\infty}\le \|v\|_{W^{1,1}}$ and by Fernique's theorem.
\end{remark}


\subsubsection{BIP for the flux function}
We now consider the inverse problem of finding the velocity function $w$ (and hence the flux function) given the initial field and finite data
\begin{equation}\label{eq:ip_flux}
    y_j = z(t_j) + \xi_j,\quad~~ j = 1, \ldots ,J,
\end{equation}
where $\xi_j\sim \mathcal{N}(0,\gamma^2)$ and are independent. Note that each $z(t_j)$ now depends implicitly on the unknown $w$. To be consistent with the notation of Theorem \ref{t:GenBay} and \ref{t:wpBay} we still set here $u=w$ and use the notation $G$ for the mapping $u\mapsto y$. We have the following well-posedness result whose proof, thanks to the continuity of the forward map $w \mapsto z$ (Theorem \ref{thm:rate_flux}), is very similar to that of Theorem \ref{t:LwpBay} and will be omitted here.

\begin{theorem}\label{t:fLwpBay}
Suppose $y\in\R^J$ is given by \eqref{eq:ip_flux}. Let $\mu_0(\Lambda)=1$ where $\Lambda$ with $Lip(\R)$ is as given in (\ref{Lip_space}). Then  the posterior measure $\mu^y$ given by 
\begin{equation*}
    \frac{d\mu^y}{d\mu_0}(u)=\frac{1}{Z}\exp\big(-\frac{1}{2\gamma^2}|y-G(u)|^2\big),
\end{equation*}
with  $Z=\int_{\Lambda} \exp\big(-\frac{1}{2\gamma^2}|y-G(u)|^2\big)\,\ud\mu_0$ is well-defined and continuous in Hellinger distance with respect to $y$.
\end{theorem}

We now show that $\mu^y$ is stable with respect to appropriate perturbations of the forward operator, provide that the shocks behave like
counterparts in traffic flow.

\begin{theorem}\label{thm:p_rate_flux}
Let the assumptions of Theorem \ref{t:fLwpBay} hold. Suppose also that for any $u\in \Lambda$, the approximating sequence $\{u^N\} \subset \Lambda$ satisfies
\begin{equation*}
    \|u^N-u\|_{\rm Lip} \le \psi(N), \quad \psi(N) \to 0 \text{ as } N\to \infty.
\end{equation*}
Assume that $\int_\Lambda \|u\|_{\rm Lip}\,\mu_0(\ud u)<\infty$ and $\mu$ and $\mu^N$ are given as in \eqref{e:posterior} and \eqref{e:Nposterior} with $G^N=G(u^N)$.

Then 
\begin{equation*}
    d_{\mathrm{Hell}}(\mu^N,\mu) \le C\sqrt{\psi(N)},
\end{equation*}
as $N\to \infty$. 
\end{theorem}

\begin{remark} A related situation to the estimation of the flux function considered here is the problem considered in \cite{HPR14} where the authors study the following inhomogeneous scalar conservation law
\begin{equation}\label{eq:inhom_tf}
    \partial_t \rho(x,t) + \partial_x (k(x)g(\rho(x,t))) = 0,
\end{equation}
with some appropriate function $k: \mathbb{R} \to \mathbb{R}$. In the context of traffic flow, the 
function $k$ represents external factors that influence the traffic flow. Such factors may be interpreted as road conditions or the presence of obstacles on the road. 
In particular, when $k$ is constant, they show that one can recover, using Tikhonov regularisation, a piecewise linear interpolation $g_\nu$ that approximates $g$ in the sense that the solutions corresponding to $g$ and $g_\nu$ are close in $L^1$. However, they assume that the whole solution is known at (almost) every point (or almost every point except some certain interval) and the initial data is chosen to be piecewise constant only.
\end{remark}

\subsection{BIP with discrete measurements of entropy solution}
\label{sec:app_Euler}

In this section, we suppose that we make a finite number of noisy observations $y=(y_1~y_2~\dots~y_J)^T$  from the solution $v$ of \eqref{scl}-\eqref{ini}. We discuss the well-posedness and also the approximation of the Bayesian inverse problem \eqref{eq:BIP} posed in Section \ref{subsec:bip} where $u$ here is  the initial field $v_0$ (the case where the flux is unknown can be treated similarly). We set $X = L^1\cap BV(\mathbb{R})$. 
The observation maps of the solution and of the approximate solution are given as
\begin{equation}\label{e:EGGN}
    G(u) = \mathcal{O}\circ S_e(u), \quad \text{and } G^N(u) = \mathcal{O}\circ S_e^N(u),
\end{equation}
where $\mathcal{O}: X \to \mathbb{R}^J$ is the evaluation operator, given by
\begin{equation*}
    \mathcal{O}(v) = \{v(x_j,t_j)\}_{j=1}^J.
\end{equation*}
 and $S_e, S_e^N: X \to X$ are the forward solution operator of \eqref{scl}-\eqref{ini} and its approximation respectively. 
This approximation may be coming from any computational method involved in solving the forward problem (including front tracking and vanishing viscosity) or a result of some disturbance in the initial condition. We show here that the Bayesian inverse problem is well-defined and continuous in data $y$ because the measurability of the forward map still holds. Lack of continuity of the forward map however leads to weaker approximation properties of the posterior compared to the case where the data came from tracking particle trajectories. We write
\begin{equation}\label{eq:EuData}
y_j=v(x_j,t_j)+\xi_j,\quad~~j=1,\dots J,
\end{equation}
with $\xi_j\sim \mathcal{N}(0,\gamma^2)$ and independent,
and let
\begin{align}\label{e:EG(u)}
    \begin{array}{l}
         u := v(\cdot,0) \\
         G(u):=(G_1(u),~\dots~, G_J(u)),\mbox{~ with ~}G_j(u)=v(x_j,t_j).
    \end{array}
\end{align}
We have the following well-posedness result.
\begin{theorem}\label{t:EwpBay}
Suppose $y\in\R^J$ is given by \eqref{eq:EuData}. Let $\mu_0(X)=1$, with $X=L^1\cap BV(\R)$. Then  the posterior measure $\mu^y$ given by 
\begin{equation*}
    \frac{d\mu^y}{d\mu_0}(u)=\frac{1}{Z}\exp\big(-\frac{1}{2\gamma^2}|y-G(u)|^2\big),
\end{equation*}
with  $Z=\int_{X} \exp\big(-\frac{1}{2\gamma^2}|y-G(u)|^2\big)\,\ud\mu_0$ is well-posed and continuous in Hellinger distance with respect to $y$.
\end{theorem}

\begin{proof}
The positivity and boundedness of $G$ imply $0<Z<1$. As the likelihood is continuous in $y$ it remains to show that $G=\mathcal{O}\circ S_e:L^1\cap BV\to\R^J$ is measurable. Since $\mathcal{O}$ is continuous, it is sufficient to show that  $S_e$ is measurable. By \cite[Theorem 6.2.3]{Daf16},  $S_e:L^1\to L^1$ is continuous. Hence for an open set $A\subset L^1$ we have that $S_e^{-1}(A)$ is open in $L^1$. Then, since $S_e:BV\to BV$ is bounded \cite[Theorem 2.14]{HR15}, we have 
$$
S_e^{-1}(A\cap BV)=S_e^{-1}(A)\cap BV\;\mbox{ is open in }\;L^1\cap BV.
$$
Noting that the Borel sigma algebra $\mathcal{B}(L^1\cap BV)=\{A\cap BV:A\in \mathcal{B}(L^1)\}$, we conclude the measurability of $S_e$ over $L^1\cap BV$. Now all conditions of Theorem \ref{t:GenBay} and \ref{t:wpBay} are satisfied and the result follows.
\end{proof}

The approximations of the posterior is much more complicated in this case. The general approximation theory of Bayesian inverse problems (\cite[Section 4.2]{DS16}) fails since the stability of the observation operator is no longer satisfied. Nevertheless, below we discuss a weaker result where one still has an error bound for the approximate posterior. 

In practical applications, one may assume that the measurement device can detect 
the shocks of strength bigger than $\varepsilon$ for some given $\varepsilon > 0$ and avoid them, where the shock strength is defined as 
\begin{equation*}
    |v(x(t)+,t) - v(x(t)-,t)|,
\end{equation*}
for some shock curve $x = x(t)$. Then the measurements may be made in the areas where no shocks or only shocks of strength smaller than $\varepsilon$ are present. We note that by \cite{BL99}, for any $\varepsilon > 0$, shocks of strength greater than $\varepsilon$ are finitely many (and hence have measure zero). Moreover, the shock strength and location of a shock of large strength do not be affected much by small changes in the initial field. In the case that the initial field and flux function are smooth ($C^k$ with $k\ge 3$), generically, solutions produce only a finite set of shocks in a given bounded domain, thanks to a regularity result discovered by Schaeffer \cite{Sch73} (see also Dafermos \cite{Daf85} for considerable improvements). Therefore, in reality, the chance of having data from certain positions with no shocks or shocks of strength smaller than some given $\varepsilon$ is likely. Let us formalise these discussions below.

Define 
$$
\mathcal{O}_\varepsilon:=\{(x,t)\in\R\times[0,T]:|v(x+,t)-v(x-,t)|>\varepsilon\},
$$
the set of shocks of strength larger than $\varepsilon$. 
For a small $\delta>0$ let the $\mathcal{O}_\varepsilon^\delta$ denote the $\delta$-neighbourhood of $\mathcal{O}_\varepsilon$.
We have the following approximation result.

\begin{theorem}
Let assumptions of Theorem \ref{t:EwpBay} hold.
Assume, for some small $\varepsilon,\delta > 0$, that the data $y=(y_1,\dots,y_J)$ is given as in \eqref{eq:EuData} with $(x_j,t_j)$, $j=1,\dots,J$, lying outside $\mathcal{O}_\varepsilon^\delta$. 
Let $\mu$ and $\mu^N$ be given by \eqref{e:posterior} and \eqref{e:Nposterior} respectively with $G$ and $G^N$ defined as in \eqref{e:EGGN} and \eqref{e:EG(u)}. 
Suppose that $S_e$ and $S_e^N$ in \eqref{e:EGGN} satisfy, for any fixed $u\in X$,  $\|S_e(u)-S_e^N(u)\|_{L^1}\to 0$ as $N\to\infty$. Then
\begin{equation*}
    d_{\rm Hell}(\mu^N, \mu) \le C\varepsilon,
\end{equation*}
for $N$ sufficiently large and with $C$ depending on $J$ and the covariance of the measurement noise.
\end{theorem}

\begin{proof}
For a given $u\in X$, by \cite[Theorem 5]{BL99}, for large enough $N$, the shocks of strength larger than $\varepsilon$ of $S_e^N(u)$ are in a $\delta$-neighbourhood of shocks of strength larger than $\varepsilon$ of $S_e(u)$. Therefore outside $\mathcal{O}_\varepsilon^\delta$ both $S_e(u)$ and $S_e^N(u)$ only have shocks of strength at most $\varepsilon$. 

For each $N$, take an integrable, bounded and continuous function, denoted by $\tilde S^N(u)$, such that $\tilde S^N(u)$ lies in the $\varepsilon$-neighbourhood of $S_e^N(u)$. Since $S_e^N(u)$ converges to $S_e(u)$ almost everywhere along a subsequence, $\tilde S^N(u)$ converges everywhere, along a subsequence, to a bounded continuous function, denoted by $\tilde S(u)$, such that $\tilde S(u)$ lies in $\varepsilon$-neighbourhood of $S_e(u)$. We can always choose $\tilde S^N(u)$ differently in the $\varepsilon$-neighbourhood of $S_e^N(u)$ such that the whole sequence $\tilde S^N(u)$ converges everywhere to $\tilde S(u)$. 
It therefore follows that
$$ \|S_e^N(u)- S_e(u)\|_{L^\infty}\le 2\varepsilon, $$
for $N$ sufficiently large. The result then follows arguing along the lines of the proof of Theorem \ref{thm:p_conv-0-rate}.
\end{proof}

In general, it may not be feasible to know a priori if shocks (of any strength) appear in certain areas, and the shock strengths may be too small that one cannot observe them. However, if the measurement devices can detect shocks of strength larger than some given $\varepsilon$, the assumption that the data $y$ is collected away from the set $\mathcal{O}_\varepsilon$ can be justified. 

\begin{remark}
If one can collect the data using the ball evaluation operator, that is, $\mathcal{O} = \{\mathcal{O}_j\}_{j=1}^J: X \to \mathbb{R}^J$, where 
\begin{equation*}
    \mathcal{O}_j(v) : = \int_{B_r(x_j)} v \ud x,
\end{equation*}
here we denote $B_r(x_j) = \{x: |x - x_j| < r\}$. Then, any small perturbation of $u$ in $X$ translates to a small perturbation of $\mathcal{O}$ in $L^\infty$, thanks to Theorem \ref{thm:Kruzkov} and Theorem \ref{thm:Lucier}. 
Therefore, the condition \eqref{eq:cont_FWmap} follows and so does Theorem \ref{thm:app_gen}. In fact, by using this method, we can estimate any unknown initial field $v_0 \in L^1 \cap L^\infty$ given a locally Lipschitz $f$, and any unknown locally Lipschitz continuous flux function for a given $v_0 \in X$. We note that the flux function need not be convex here (see Theorem \ref{thm:Lucier}). 
\end{remark}


\section*{Acknowledgments} 
The authors are grateful to Konstantinos Koumatos and Aretha Teckentrup for helpful discussions.
MD was partially supported by a fellowship from
the Simons Foundation. The work of DLD was partly supported by the European Union’s Horizon 2020 research and innovation programme under the Marie Skłodowska-Curie grant agreement (No 642768), and the Academy of Finland (No 345720). The second author would like to thank Andrew Duncan for his encouragement and helpful discussions at the beginning of this project.

\bibliographystyle{alpha}
\bibliography{bibliography}
\end{document}